\documentclass{amsart}

\usepackage{amssymb,amsthm,amsmath,amsxtra}
\usepackage[all]{xy}
\usepackage{verbatim}

\numberwithin{equation}{section}

\theoremstyle{plain}
\newtheorem{thm}[equation]{Theorem}
\newtheorem{prop}[equation]{Proposition}
\newtheorem{lem}[equation]{Lemma}
\newtheorem{cor}[equation]{Corollary}
\newtheorem{prob}[equation]{Problem}
\newtheorem{conj}[equation]{Conjecture}
\newtheorem*{prob*}{Problem}
\newtheorem*{thm*}{Theorem}

\theoremstyle{definition}
\newtheorem{alg}[equation]{Algorithm}
\newtheorem{defn}[equation]{Definition}

\theoremstyle{remark}
\newtheorem{exm}[equation]{Example}
\newtheorem{rmk}[equation]{Remark}
\newtheorem{ques}[equation]{Question}

\newenvironment{enumroman}
{\begin{enumerate}}
{\end{enumerate}}

\newenvironment{enumalg}
{\begin{enumerate}}
{\end{enumerate}}

\newenvironment{enumalgalph}
{\begin{enumerate}}
{\end{enumerate}}

\entrymodifiers={+!!<0pt,\fontdimen22\textfont2>}

\DeclareMathOperator{\opchar}{char}

\DeclareMathOperator{\disc}{disc}
\DeclareMathOperator{\End}{End}

\DeclareMathOperator{\GL}{GL}

\DeclareMathOperator{\M}{M}
\DeclareMathOperator{\N}{N}
\DeclareMathOperator{\nrd}{nrd}
\DeclareMathOperator{\ord}{ord}
\DeclareMathOperator{\PGL}{PGL}
\DeclareMathOperator{\rad}{rad}

\DeclareMathOperator{\sqrad}{sqrad}
\DeclareMathOperator{\trd}{trd}
\DeclareMathOperator{\Tr}{Tr}

\newcommand{\C}{\mathbb C}
\newcommand{\F}{\mathbb F}
\newcommand{\HH}{\mathbb H}
\newcommand{\PP}{\mathbb P}
\newcommand{\Q}{\mathbb Q}
\newcommand{\R}{\mathbb R}
\newcommand{\Z}{\mathbb Z}

\newcommand{\calO}{\mathcal O}

\newcommand{\fraka}{\mathfrak{a}}
\newcommand{\frakb}{\mathfrak{b}}
\newcommand{\frakc}{\mathfrak{c}}
\newcommand{\frakd}{\mathfrak{d}}
\newcommand{\frakD}{\mathfrak{D}}
\newcommand{\frakg}{\mathfrak{g}}

\newcommand{\frakp}{\mathfrak{p}}
\newcommand{\frakq}{\mathfrak{q}}

\newcommand{\psmod}[1]{~(\textup{\text{mod}}~{#1})}

\newcommand{\legen}[2]{\left(\frac{#1}{#2}\right)}

\newcommand{\quat}[2]{\displaystyle{\biggl(\frac{#1}{#2}\biggr)}}
\newcommand{\issq}[2]{\displaystyle{\biggl\{\frac{#1}{#2}\biggr\}}}

\newcommand{\probsf}[1]{\textup{(\textsf{#1})}}

\newcommand{\la}{\langle}
\newcommand{\ra}{\rangle}

\newcommand{\Lambdabar}{\overline{\Lambda}}

\begin{document}

\title[Identifying the matrix ring]{Identifying the matrix ring: algorithms for quaternion algebras and quadratic forms}

\author{John Voight}
\address{Department of Mathematics and Statistics, University of Vermont, 16 Colchester Ave, Burlington, VT 05401, USA}
\email{jvoight@gmail.com}
\date{\today}



\begin{abstract}
We discuss the relationship between quaternion algebras and quadratic forms with a focus on computational aspects.  Our basic motivating problem is to determine if a given algebra of rank $4$ over a commutative ring $R$ embeds in the $2 \times 2$-matrix ring $\M_2(R)$ and, if so, to compute such an embedding.  We discuss many variants of this problem, including algorithmic recognition of quaternion algebras among algebras of rank $4$, computation of the Hilbert symbol, and computation of maximal orders.
\end{abstract}

\maketitle
\tableofcontents

Since the discovery of the division ring of quaternions over the real numbers by Hamilton, and continuing with work of Albert and many others, a deep link has been forged between quadratic forms in three and four variables over a field $F$ and quaternion algebras over $F$.  Starting with a \emph{quaternion algebra} over $F$, a central simple $F$-algebra of dimension $4$, one obtains a quadratic form via the reduced norm (restricted to the trace zero subspace); the split quaternion algebra over $F$, the $2 \times 2$-matrix ring $\M_2(F)$, corresponds to an isotropic quadratic form, one that represents zero nontrivially.  (Conversely, one recovers the quaternion algebra via the Clifford algebra of the quadratic form.)  In this article, we give an exposition of this link relating quaternion algebras and quadratic forms from an explicit, algorithmic perspective and in a wider context.

Let $R$ be a noetherian, commutative domain.  
We say that $R$ is \emph{computable} if there exists an encoding of $R$ into bits with algorithms to perform ring operations in $R$ and to test if an element of $R$ is zero.  The following basic algorithmic problem, along with its many variants, forms the core of this article.  (See \S 1 for further definitions and algorithmic specifications.)

\begin{prob*}[\textsf{IsMatrixRing}]
Given a computable domain $R$ and an $R$-algebra $\calO$ of rank $4$, determine if $\calO$ embeds in $\M_2(R)$ and, if so, compute an explicit embedding $\calO \hookrightarrow \M_2(R)$ of $R$-algebras.
\end{prob*}

The problem (\textsf{IsMatrixRing}) captures in an important way the link between quadratic forms and quaternion algebras.  In the simplest case where $R=F$ is a field---when such an embedding is necessarily an isomorphism---this problem corresponds to asking if a ternary quadratic form over $F$ represents zero nontrivially, and for this reason it arises in a wide variety of situations.  When $F$ is a local field, this problem corresponds to the computation of the Hilbert symbol.  In the case where $R$ is a local ring, it corresponds to the computation of an (explicit) integral splitting of a quaternion order and thereby appears as a foundational step in many algorithms in arithmetic geometry (as in work of Kirschmer and the author \cite{KirschmerVoight}).  Finally, when $R$ is a Dedekind domain, roughly speaking, the problem of approximating (\textsf{IsMatrixRing}) naturally gives rise to the problem of computing a maximal order containing $\calO$.  In these and other ways, therefore, the problem (\textsf{IsMatrixRing}) will serve as kind of unifying and motivating question.

In \S 1, we introduce the basic terminology we will use throughout concerning computable rings and quaternion algebras.  In \S 2, we consider algebras equipped with a standard involution and we exhibit an algorithm to test if an $F$-algebra $B$ has a standard involution.  In \S 3, we relate algebras with a standard involution to quadratic forms via the reduced norm; we introduce the theory of quadratic forms over local PIDs, providing an algorithm to compute a \emph{normalization} of such a form.  As a consequence, we exhibit an algorithm to test if an $F$-algebra $B$ is a quaternion algebra and, if so, to compute \emph{standard generators} for $B$.  With these reductions, we turn in \S 4 to Problem (\textsf{IsMatrixRing}) for quaternion algebras and prove that this problem is deterministic polynomial-time equivalent to the problem of determining if a conic defined over $F$ has an $F$-rational point (and, if so, to exhibit one).   

In \S 5, we consider Problem (\textsf{IsMatrixRing}) in the case where $F$ is a local field, which corresponds to the computation of the Hilbert symbol; in \S 6 we treat the more delicate case of a local dyadic field, and putting these together prove that there is a deterministic polynomial-time algorithm to compute the Hilbert symbol (Theorem \ref{evenhilb}).  We thereby exhibit an algorithm to compute the generalized Jacobi symbol for computable Euclidean domains.  In \S 7, we turn to the case of a Dedekind domain $R$ and relate Problem (\textsf{IsMatrixRing}) to the problem of computing a maximal $R$-order; we prove that the problem of computing a maximal order for a quaternion algebra $B$ over a number field $F$ is probabilistic polynomial-time equivalent to the problem of factoring integers.  Finally, in \S 8, we consider the problem (\textsf{IsMatrixRing}) over $\Q$, and show that recognizing the matrix ring is  deterministic polynomial-time equivalent to the problem of quadratic residuosity.

Many of the results in this paper fit into the more general setting of semisimple algebras; however, we believe that the special link to quadratic forms, along with the wide application of quaternion algebras (analogous to that of quadratic field extensions), justifies the specialized treatment they are afforded here.

The author would like to thank his Ph.D.\ advisor Hendrik Lenstra for his many helpful comments, the \textsf{Magma} group at the University of Sydney for their support while writing this paper, and David Kohel for his valuable input.  We are indebted to Carl Pomerance for the citation \cite{AdPomRum} and would like to thank Asher Auel, Jonathan Hanke, Kate Thompson, and the referee for helpful corrections and suggestions.  Some of the results herein occur in the author's Ph.D. thesis \cite{Voight}.  Writing this paper was partially supported by the National Security Agency under Grant Number H98230-09-1-0037 and the National Science Foundation under Grant No.\ DMS-0901971.

\section{Rings and algebras}

We begin by introducing some notation and background that will be used throughout.  Let $R$ be a commutative, noetherian domain (with $1$), and let $F$ be the field of fractions of $R$.  

Let $\calO$ be an \emph{$R$-algebra}, an associative ring with $1$ equipped with an embedding $R \hookrightarrow \calO$ of rings (taking $1 \in R$ to $1 \in \calO$) whose image lies in the center of $\calO$; we identify $R$ with its image under this embedding.  We will assume without further mention that $\calO$ is a finitely generated, projective (equivalently, locally free) $R$-module of rank $n \in \Z_{\geq 1}$.

\subsection*{Computable rings and algebras}

We will follow the conventions of Lenstra \cite{Lenstra} for rings and algorithms, with the notable exception that we do not require all rings to be commutative.

A domain $R$ is \emph{computable} if $R$ comes equipped with a way of encoding elements of $R$ in bits (i.e.\ the elements of $R$ are recursively enumerable, allowing repetitions) along with deterministic algorithms to perform ring operations in $R$ (addition, subtraction, and multiplication) and to test if $x=0 \in R$; a ring is \emph{polynomial-time computable} if these algorithms run in polynomial time (in the bit size of the input).  A field is \emph{computable} if it is a computable ring and furthermore there exists an algorithm to divide by a nonzero element.  For precise definitions and a thorough survey of the subject of computable rings we refer to Stoltenberg-Hansen and Tucker \cite{SHT} and the references contained therein.

\begin{exm}
A domain $R$ which is the localization of a ring which is finitely generated over its prime ring is computable by the theory of Gr\"obner bases \cite{GathenGerhard}.  For example, any finitely generated algebra over $\Z$ or $\Q$ (without zerodivisors, since we restrict to domains) is computable, and in particular the coordinate ring of any integral affine variety over a finitely generated field is computable.
\end{exm}

\begin{exm}
If $R$ is a computable domain, then $F$ is a computable field if elements are represented in bits as pairs of elements of $R$ in the usual way.  
\end{exm}

\begin{rmk}
Inexact fields (e.g.\ local fields, such as $\Q_p$ or $\R$) are not computable, since they are uncountable!  However, see the discussion in \S \ref{localfields} for the use of a computable subring which works well in our situation.
\end{rmk}

\begin{exm}
A number field $F$ is computable, specified by the data of the minimal polynomial of a primitive element (itself described by the sequence of its coefficients, given as rational numbers); elements of $F$ are described by their standard representation in the basis of powers of the primitive element \cite[\S 4.2.2]{Cohen}.  For a detailed exposition of algorithms for computing with a number field $F$, see Cohen \cite{Cohen,Cohen2} and Pohst and Zassenhaus \cite{PZ}.  
\end{exm}

\begin{rmk} \label{globalfunfields}
\emph{Global function fields}, i.e.\ finite extensions of $k(T)$ with $k$ a finite field, can be treated in a parallel fashion to number fields.  Unfortunately, at the present time the literature is much less complete in providing a suite of algorithms for computing with integral structures in such fields---particularly in the situation where one works in a relative extension of such fields---despite the fact that some of these algorithms have already been implemented in \textsf{Magma} \cite{Magma} by Hess \cite{Hess}.  Therefore, in this article we will often consider just the case of number fields and content ourselves to notice that the algorithms we provide will generalize with appropriate modifications to the global function field setting.
\end{rmk}

Throughout this article, when discussing algorithms, we will assume that the domain $R$ and its field of fractions $F$ are computable.

Let $B$ be a $F$-algebra with $\dim_F B=n$ and basis $e_1,e_2,\dots,e_n$ (as an $F$-vector space), and suppose $e_1=1$.  A \emph{multiplication table} for $B$ is a system of $n^3$ elements $(c_{ijk})_{i,j,k=1,\dots,n}$ of $F$, called \emph{structure constants}, such that multiplication in $B$ is given by
\[ e_i e_j = \sum_{k=1}^n c_{ijk} e_k \]
for $i,j \in \{1,\dots,n\}$.  

An $F$-algebra $B$ is represented in bits by a multiplication table and elements of $F$ are represented in the basis $e_i$.  Note that basis elements in $B$ can be multiplied directly by the multiplication table but multiplication of arbitrary elements in $B$ requires $O(n^3)$ arithmetic operations (additions and multiplications) in $F$; in either case, note the output is of polynomial size in the input for fixed $B$.

\begin{rmk}
We have assumed that $B$ is associative as an $F$-algebra; however, this property can be verified by simply checking the associative law on a basis.
\end{rmk}

\begin{rmk}
We require that the element $1$ be included as a generator of $B$, since by our definition an $F$-algebra is equipped with an embedding $F \hookrightarrow B$.  This is not a serious restriction, for the equations which uniquely define the element $1$ in $B$ are linear equations and so $1 \in B$ can be (uniquely) recovered by linear algebra over $F$.  (And an algebra without $1$ embeds inside an algebra with $1$.)
\end{rmk}


An $R$-algebra $\calO$ is represented in bits by the $F$-algebra $B=\calO \otimes_R F$ and a set of $R$-module generators $x_1,\dots,x_m \in B$ with $x_1=1$.  A morphism between $R$-algebras is represented by the underlying $R$-linear map, specified by a matrix in the given sets of generators for the source and target.  

\subsection*{Quaternion algebras}

We refer to Vign\'eras \cite{Vigneras} and Reiner \cite{Reiner} for background relevant to this section.  

An $F$-algebra $B$ is \emph{central} if the center of $B$ is equal to $F$, and $B$ is \emph{simple} if the only two-sided ideals of $B$ are $(0)$ and $B$ (or equivalently that any $F$-algebra homomorphism with domain $B$ is either the zero map or injective).

\begin{rmk} \label{computecenter}
One can compute the center of $B$ by solving the $n$ linear equations $xe_i=e_ix$ for $x=x_1e_1 + \dots + x_ne_n$ and thereby, for example, verify that $B$ is central.
\end{rmk}

\begin{defn}
A \emph{quaternion algebra} $B$ over $F$ is a central simple $F$-algebra with $\dim_F B=4$.  
\end{defn}

An $F$-algebra $B$ is a quaternion algebra if and only if there exist $i,j \in B$ which generate $B$ as an $F$-algebra such that
\begin{equation} \label{stdgennotchar2}
i^2=a, \quad j^2=b, \quad ji=-ij
\end{equation}
with $a,b \in F^\times$ if $\opchar F \neq 2$, and
\begin{equation} \label{stdgenchar2}
i^2+i=a,\quad j^2=b,\quad ji=(i+1)j
\end{equation}
with $a \in F$ and $b \in F^\times$ if $\opchar F = 2$.  We give an algorithmic proof of this equivalence in \S 3.  We accordingly denote an algebra (\ref{stdgennotchar2})--(\ref{stdgenchar2}) by $B=\quat{a,b}{F}$, say that $B$ is in \emph{standard form}, and call the elements $i,j$ \emph{standard generators}.  Note that $B$ has basis $1,i,j,ij$ as an $F$-vector space, so indeed $\dim_F B=4$.

\begin{exm} \label{M2Fisom11}
The ring $\M_2(F)$ of $2 \times 2$-matrices with coefficients in $F$ is a quaternion algebra over $F$.  Indeed, we have $\quat{1,1}{F} \cong \M_2(F)$ with $j \mapsto \begin{pmatrix} 0 & 1 \\ 1 & 0 \end{pmatrix}$ and
\[ i \mapsto \begin{pmatrix} 1 & 0 \\ 0 & -1 \end{pmatrix}\quad \text{or}\quad i \mapsto \begin{pmatrix} 0 & 1 \\ 1 & 1 \end{pmatrix} \] 
according as $\opchar F \neq 2$ or $\opchar F = 2$. 

Every quaternion algebra over a separably (or algebraically) closed field $F$ is isomorphic to $\M_2(F)$.  
\end{exm}

\begin{exm}
The $\R$-algebra $\HH$, generated by $i,j$ satisfying $i^2=j^2=(ij)^2=-1$ is the usual division ring of quaternions over $\R$.  Every quaternion algebra over $\R$ is isomorphic to either $\M_2(\R)$ or $\HH$, according to the theorem of Frobenius.
\end{exm}

Let $B$ be an $F$-algebra.  An \emph{$R$-order} in $B$ is a subring $\calO \subset B$ that is finitely generated as an $R$-module and such that $\calO F = B$.  We see that an $R$-algebra $\calO$ is an $R$-order in $B=\calO \otimes_R F$, and we will use this equivalence throughout, sometimes thinking of $\calO$ as an $R$-algebra on its own terms and at other times thinking of $\calO$ as arising as an order inside an algebra over a field.

A \emph{quaternion order} over $R$ is an $R$-order in a quaternion algebra $B$ over $F$.  Equivalently, an $R$-algebra $\calO$ is a quaternion order if $B=\calO \otimes_R F$ is a quaternion algebra over $F$.  

\begin{exm}
$\M_2(R)$ is a quaternion order in $\M_2(F)$.  

If $a,b \in R$ then $\calO=R \oplus Ri \oplus Rj \oplus Rij$ is a quaternion order in $B=\quat{a,b}{F}$.  So for example $\Z \oplus \Z i \oplus \Z j \oplus \Z ij$ is a $\Z$-order in the rational Hamiltonians $B=\quat{-1,-1}{\Q}$.
\end{exm}

Further examples of quaternion orders will be defined in the next section (see Lemma \ref{SibR}).

\subsection*{Modules over Dedekind domains}

Let $R$ be a Dedekind domain, an integrally closed (noetherian) domain in which every nonzero prime ideal is maximal.  Every field is a Dedekind domain (vacuously), as is the integral closure of $\Z$ or $\F_p[T]$ in a finite (separable) extension of $\Q$ or $\F_p(T)$, respectively.  The localization of a Dedekind domain at a multiplicative subset is again a Dedekind domain.  If $R$ is the ring of integers of a number field, then we call $R$ a \emph{number ring}.

Over a Dedekind domain $R$, every projective $R$-module $M$ can be represented as the direct sum of projective $R$-modules of rank $1$, which is to say that there exist projective (equivalently, locally principal) $R$-modules $\fraka_1,\dots,\fraka_n \subset F$ (also known as \emph{fractional ideals} of $R$) and elements $x_1,\dots,x_n \in M$ with $\fraka_1=R$ and $x_1=1$ such that
\[ M = \fraka_1 x_1 \oplus \dots \oplus \fraka_n x_n; \]
we say then that the elements $x_i$ are a \emph{pseudobasis} for $M$ with \emph{coefficient ideals} $\fraka_i$.  More generally, if $M = \fraka_1 x_1 + \dots + \fraka_m x_m$ (the sum not necessarily direct), then we say the elements $x_i$ are a \emph{pseudogenerating set} for $M$ (with \emph{coefficient ideals} $\fraka_i$).  

In fact, the above characterization can be made computable as follows.

\begin{prop} \label{HMFalg}
Let $R$ be a number ring.  Then there exists an algorithm which, given a projective $R$-module $M$ specified by a pseudogenerating set, returns a pseudobasis for $M$.
\end{prop}

The algorithm in Proposition \ref{HMFalg} is a generalization of the Hermite normal form (HMF) for matrices over $\Z$; see Cohen \cite[Chapter 1]{Cohen2}.  Therefore, from now on we represent a quaternion order $\calO$ over a number ring $R$ by a pseudobasis; in such a situation, we may and do assume that $x_1=1$ (by employing the HMF).

\begin{rmk}
Recalling Remark \ref{globalfunfields}, in particular there seems to be no comprehensive reference for results akin to Proposition \ref{HMFalg} in the global function field case.
\end{rmk}

\section{Standard involutions and degree} \label{secdegree}

Quaternion algebras, or more generally algebras which have a standard involution, possess a quadratic form called the reduced norm.  In this section, we discuss this association and we give an algorithm which verifies that an algebra has a standard involution.  As a reference, see Jacobson \cite[\S 1.6]{Jacobson}, Knus \cite{Knus}, and work of the author \cite{Voightlowrank}.

In this section, let $R$ be an integrally closed (noetherian) domain with field of fractions $F$.  Let $\calO$ be an $R$-algebra and let $B=\calO \otimes_R F$.

\subsection*{Degree}

We first generalize the notion of degree from field extensions to $R$-algebras.

\begin{defn}
The \emph{degree} of $x \in \calO$ over $R$, denoted $\deg_R(x)$, is the smallest positive integer $n$ such that $x$ satisfies a monic polynomial of degree $n$ with coefficients in $R$.  The \emph{degree} of $\calO$ over $R$, denoted $\deg_R(\calO)$, is the smallest positive integer $n$ such that every element of $\calO$ has degree at most $n$.
\end{defn}

Every $x \in \calO$ satisfies the characteristic polynomial of (left) multiplication by $x$ on a set of generators for $\calO$ as an $R$-module, and consequently $\deg_R(\calO) < \infty$ (under our continuing hypothesis that $\calO$ is projective of finite rank).  

\begin{lem} \label{degROdegFB}
We have $\deg_R(\calO)=\deg_F(B)$.
\end{lem}

\begin{proof}
Since $\calO$ is finitely generated as an $R$-module and $R$ is noetherian, the $R$-submodule $R[x] \subset \calO$ is finitely generated, so $x$ is integral over $R$.  Since $R$ is integrally closed, the minimal polynomial of $x \in \calO$ over $F$ has coefficients in $R$ by Gauss' lemma, so $\deg_R(x)=\deg_F(x)$ and thus $\deg_R(\calO) \leq \deg_F(B)$.  On the other hand, if $y \in B$ then there exists $0 \neq d \in R$ such that $x=yd \in \calO$ so $\deg_F(x)=\deg_F(y)=\deg_R(y)$ so $\deg_F(B) \leq \deg_R(\calO)$.
\end{proof}

From the lemma, we need only consider the degree of an algebra over a field.

\begin{exm} \label{deg2exm}
$B$ has degree $1$ if and only if $B=F$.  

If $K$ is a separable field extension of $F$ with $\dim_F K=n$, then $K$ has degree $n$ as a $F$-algebra (in the above sense) by the primitive element theorem.  

If $\dim_F B=n$, then $B$ has degree at most $n$ but even if $B$ is commutative one may still have $\deg_F(B) < \dim_F B$: for example, $B=F[x,y,z]/(x,y,z)^2$ has rank $4$ over the field $F$ but has degree $2$.
\end{exm}

\subsection*{Standard involutions}

We will see in a moment that quaternion orders and algebras are algebras of degree $2$; this will be a consequence of the fact that they possess a standard involution.  Indeed, the link between algebras with an involution and quadratic forms forms the heart of much important work \cite{bookofinvolutions}.

\begin{defn}
An \emph{anti-automorphism} of $\calO$ is an $R$-linear map $\overline{\phantom{x}}:\calO \to \calO$ with $\overline{1}=1$ and $\overline{xy}=\overline{y}\,\overline{x}$ for all $x \in \calO$.  An \emph{involution} is an anti-automorphism such that $\overline{\overline{x}}=x$ for all $x \in \calO$.  An involution is \emph{standard} if $x\overline{x} \in R$ for all $x \in \calO$.
\end{defn}  

Note that if $x\overline{x} \in R$ for all $x \in \calO$, then $(x+1)(\overline{x}+1)=x\overline{x}+(x+\overline{x})+1 \in R$ and hence $x+\overline{x} \in R$ for all $x \in \calO$ as well.  Note that $\overline{x}x=x\overline{x}$ for all $x \in \calO$ since $x(x+\overline{x})=(x+\overline{x})x$ (and $R$ is central in $\calO$).

\begin{exm}
If $\calO=\M_n(R)$, then the transpose map is an anti-automorphism which is standard if and only if $n=1$; the adjoint map is a standard involution for $n \leq 2$ but is not $R$-linear for $n \geq 3$.
\end{exm}

Suppose now that $\calO$ has a standard involution $\overline{\phantom{x}}$.  Then we define the \emph{reduced trace} and \emph{reduced norm}, respectively, to be the maps
\begin{align*}
\trd:\calO &\to R & \nrd:\calO &\to R \\
x &\mapsto x+\overline{x} & x &\mapsto x\overline{x}=\overline{x}x
\end{align*}

We have
\begin{equation} \label{quadeqnsatisfied}
x^2-\trd(x)x+\nrd(x)=x^2-(x+\overline{x})x+x\overline{x}=0
\end{equation}
for all $x \in \calO$.  It follows that if $\calO$ has a standard involution then either $\calO=R$ (so the standard involution is the identity and $\calO=R$ has degree $1$) or $\calO$ has degree $2$.  

\begin{exm} \label{quatstdinvo}
Let $B=\quat{a,b}{F}$ be a quaternion algebra over $F$.  Then $B$ has a standard involution, defined as follows.  For $x=t+ui+vj+wk$, we have
\[ \overline{x}=t-ui-vj-wk \]
so $\trd(x)=2t$ and $\nrd(x)=t^2-au^2-bv^2+abw^2$ if $\opchar F \neq 2$ and
\[ \overline{x}=t+(u+1)i+vj+wk \]
so $\trd(x)=2u$ and $\nrd(x)=t^2+tu+au^2 + bv^2+bvw+abw^2$ if $\opchar F = 2$.  
\end{exm}

\begin{lem} \label{BOFstdinv}
$\calO$ has a standard involution if and only if $B=\calO \otimes_R F$ has a standard involution.
\end{lem}

\begin{proof}
If $\calO$ has a standard involution, we obtain one on $B$ by extending $F$-linearly.  Conversely, suppose $B$ has a standard involution and let $x \in \calO$.  Then as in the proof of Lemma \ref{degROdegFB}, $x$ is integral over $R$ so its minimal polynomial over $F$ has coefficients in $R$.  If $x \in R$, then $\overline{x}=x$ and there is nothing to prove.  If $x \not \in R$, this minimal polynomial must be given by (\ref{quadeqnsatisfied}), so $\trd(x)=x+\overline{x} \in R$ and thus $\overline{x}=\trd(x)-x \in \calO$ has $x\overline{x}=\nrd(x) \in R$ as well.
\end{proof}

An $R$-algebra $S$ is \emph{quadratic} if $S$ has rank $2$ as an $R$-module.

\begin{lem} \label{quadcomm}
Let $S$ be a quadratic $R$-algebra.  Then $S$ is commutative and has a unique standard involution.  
\end{lem}

\begin{proof}
By Lemma \ref{BOFstdinv}, it suffices to prove the lemma for $K=S \otimes_R F$.  But then for any $x \in K \setminus F$ we have $K=F \oplus F x$ so $K$ is commutative.  Moreover, we have $x^2-tx+n=0$ for some unique $t,n \in F$ and so the (necessarily unique) standard involution is given by $x \mapsto t-x$, extending by $F$-linearity.  
(See also Scharlau \cite[\S 8.11]{Scharlau} for a proof of this lemma.)
\end{proof}  

\begin{cor}
If $\calO$ has a standard involution, then this involution is unique.
\end{cor}

This corollary follows immediately from Lemma \ref{quadcomm} by restricting to quadratic subalgebras $K$ of $B$.

\subsection*{Quaternion orders}

Having identified the standard involution on a quadratic algebra, we now generalize the construction of quaternion algebras (\ref{stdgennotchar2})--(\ref{stdgenchar2}) to quaternion orders.  Let $S$ be a quadratic $R$-algebra, and suppose $S$ is \emph{separable}, so the minimal polynomial of every $x \in S$ has distinct roots over the algebraic closure $\overline{F}$ of $F$.  Let $J \subset S$ be an invertible $S$-ideal (equivalently, a locally principal $S$-module) and let $b \in R \setminus \{0\}$.  We denote by $\quat{S,J,b}{R}$ the $R$-algebra $S \oplus J j$ subject to the relations $j^2=b$ and $j i=\overline{i} j$ for all $i \in S$, where $\overline{\phantom{x}}$ denotes the unique standard involution on $S$ obtained from Lemma \ref{quadcomm}.  We say that such an algebra is in \emph{standard form}.

\begin{lem} \label{SibR}
The $R$-algebra $\calO=\quat{S,J,b}{R}$ is a quaternion order.
\end{lem}

\begin{proof}
We consider $B = \calO \otimes_R F$.  Let $K = S \otimes_R F$ and let $i \in K \setminus F$.  Since $K$ is separable, if $\opchar F \neq 2$ by completing the square we may assume $i^2=a$ with $a \in F^\times$; if $\opchar F = 2$, we may assume $i^2+i=a$ with $a \in F$.  Now since $J$ is projective we have $J \otimes_R F = J \otimes_S K \cong K$ so $B \cong K \oplus Kj$ as an $F$-algebra.  Finally, since $ji=\overline{i}j=(\trd(i)-i)j$ and $\trd(i)=0,1$ according as $\opchar F \neq 2$ or not, we have identified $B$ as isomorphic to the quaternion algebra $\quat{a,b}{F}$.
\end{proof}


\subsection*{Algorithmically identifying a standard involution}

We conclude this section with an algorithm to test if an $F$-algebra $B$ (of dimension $n$) has a standard involution.

First, we note that if $B$ has a standard involution $\overline{\phantom{x}}:B \to B$, then this involution and hence also the reduced trace and norm can be computed efficiently.  Indeed, let $\{e_i\}_i$ be a basis for $B$; then $\trd(e_i) \in F$ is simply the coefficient of $e_i$ in $e_i^2$, and so $\overline{e_i}=\trd(e_i)-e_i$ for each $i$ can be precomputed for $B$; one recovers the involution on $B$ (and hence also the trace) for an arbitrary element of $B$ by $F$-linearity.  Therefore the involution and the reduced trace can be computed using $O(n)$ arithmetic operations in $F$ (with output linear in the input for fixed $B$) and the reduced norm using $O(n^2)$ operations in $F$ (with output quadratic in the input).

\begin{alg} \label{deg2}
Let $B$ be an $F$-algebra given by a multiplication table in the basis $e_1,\dots,e_n$ with $e_1=1$.  This algorithm returns \textsf{true} if and only if $B$ has a standard involution.
\begin{enumalg}
\item For $i=2,\dots,n$, let $t_i \in F$ be the coefficient of $e_i$ in $e_i^2$, and let $n_i=e_i^2-t_ie_i$.  If some $n_i \not\in F$, return \textsf{false}.
\item For $i=2,\dots,n$ and $j=i+1,\dots,n$, let $n_{ij} = (e_i+e_j)^2 - (t_i+t_j)(e_i+e_j)$.  If some $n_{ij} \not\in F$, return \textsf{false}.  Otherwise, return \textsf{true}.
\end{enumalg}
\end{alg}

\begin{proof}[Proof of correctness]
Let $F[x]=F[x_1,\dots,x_n]$ be the polynomial ring over $F$ in $n$ variables, and let $B_{F[x]} = B \otimes_F F[x]$.  Let $\xi=x_1+x_2e_2+\dots+x_ne_n \in B_{F[x]}$, and define
\[ t_\xi = \sum_{i=1}^n t_i x_i \] 
and 
\[ n_\xi = \sum_{i=1}^n n_i x_i^2 + \sum_{1\leq i<j \leq n} (n_{ij}-n_i-n_j)x_ix_j. \]
Let
\[ \xi^2-t_\xi \xi+n_\xi= \sum_{i=1}^n c_i(x_1,\dots,x_n)e_i \]
with $c_i(x) \in F[x]$.  Each $c_i(x)$ is a homogeneous polynomial of degree $2$.  The algorithm then verifies that $c_i(x)=0$ for $x \in \{e_i\}_i \cup \{e_i+e_j\}_{i,j}$, and this implies that each $c_i(x)$ vanishes identically.  Therefore, the specialization of the map $\xi \mapsto \overline{\xi}=t_\xi - \xi$ is the unique standard involution on $B$.  
\end{proof}

\begin{rmk}
Algorithm \/\textup{\ref{deg2}} requires $O(n)$ arithmetic operations in $F$, since $e_i^2$ can be computed directly from the multiplication table and hence $(e_i+e_j)^2=e_i^2+e_ie_j+e_je_i+e_j^2$ can be computed using $O(4n)=O(n)$ operations.
\end{rmk}

\begin{rmk}
Using the notation of the proof of correctness for Algorithm \ref{deg2}, it is clear that $\deg(B)=\deg(\xi)$, i.e., $\deg(B)$ is equal to the degree of the minimal polynomial of $\xi$, which can be computed as the rank of the matrix over $F[x]$ whose columns are $1,\xi,\dots,\xi^n$ using linear algebra over the field $F(x_1,\dots,x_n)$.
\end{rmk}

\section{Algebras with a standard involution and quadratic forms} \label{ssquadform}

In this section, we describe a relationship between $R$-algebras with a standard involution and quadratic forms over $R$.  The main result of this section is an algorithm which verifies that an $R$-algebra $\calO$ over a local PID is a quaternion order and, if so, exhibits standard generators for $\calO$.  Specializing, we will thereby recognize quaternion algebras over a field $F$.  We then extend this to recognizing quaternion orders over a number ring $R$.  Over fields, a reference for this section is Lam \cite{Lam}, and for more about algebras equipped with a quadratic norm form, we refer the reader to Knus \cite{Knus}.

\subsection*{Quadratic forms over rings}

We begin by defining quadratic forms over a (noetherian) domain $R$.

\begin{defn}
A \emph{quadratic form} over $R$ is a map $Q:M \to R$, where $M$ is a finitely generated projective $R$-module, such that:
\begin{enumroman}
\item $Q(ax)=a^2Q(x)$ for all $a \in R$ and $x \in M$; and
\item The map $T:M \times M \to R$ defined by
\[ T(x,y) = Q(x+y)-Q(x)-Q(y) \]
is $R$-bilinear.
\end{enumroman}
\end{defn}


A symmetric bilinear form $T:M \times M \to R$ is \emph{even} if $T(x,x) \in 2R$ for all $x \in M$.  If $T$ arises from a quadratic form, then $T$ is even, and conversely if $T$ is even and $2$ is a nonzerodivisor in $R$ then one recovers the quadratic form as $Q(x)=T(x,x)/2$.  

Let $Q:M \to R$ be a quadratic form and suppose that $M$ is free over $R$ with basis $e_1,\dots,e_n$.  The \emph{Gram matrix} of $Q$ with respect to the basis $e_1,\dots,e_n$ is the matrix $A=(T(e_i,e_j))_{i,j=1,\dots,n} \in \M_n(R)$.  The matrix $A$ has the property that $x^t A y = T(x,y)$, where we identify $x=x_1e_1 + \dots + x_ne_n$ with the column vector $(x_1,\dots,x_n)^t$, and similarly for $y$.  In particular we have $x^t A x = 2Q(x)$.  

Let $Q:M \to R$ be a quadratic form.  We say $x,y \in M$ are \emph{orthogonal} (with respect to $Q$) if $T(x,y)=0$.  

\begin{exm}
Let $\calO$ be an $R$-algebra with a standard involution $\overline{\phantom{x}}$.  Then the reduced norm $\nrd:\calO \to R$ (defined by $x \mapsto x\overline{x}$ for $x \in \calO$) is a quadratic form on $\calO$ with associated bilinear form
\begin{equation} \label{Bfm}
T(x,y) = x\overline{y}+y\overline{x}= \trd(x\overline{y})=\trd(x)y+\trd(y)x-(xy+yx)=\trd(\overline{x}y)
\end{equation}
for $x,y \in \calO$.  In particular $T(1,x)=T(x,1)=\trd(x)$.  Note that $x,y \in \calO$ are orthogonal if and only if $x\overline{y}=-y\overline{x}$, and if further $\trd(x)=\trd(y)=0$ then $\overline{x}=-x$ and $\overline{y}=-y$ so $x,y$ are orthogonal if and only if $xy=-yx$.
\end{exm}

\begin{exm}
Let $\calO_0=\{x \in \calO : \trd(x)=0\}$ be the $R$-submodule of elements of reduced trace zero.  Then $\calO/\calO_0$ is torsion-free, since if $rx \in \calO_0$ then $\trd(rx)=r\trd(x)=0$ so $\trd(x)=0$ so $x \in \calO_0$.  Thus $\calO_0$ is a projective $R$-submodule of $\calO$ and $\calO \supset R \oplus \calO_0$.  We therefore obtain a quadratic form $\nrd_0=\nrd|_{\calO_0}:\calO_0 \to R$.  
\end{exm}

If $Q:M \to R$ and $Q':M' \to R$ are quadratic forms, we define the form $Q \perp Q'$ on $M \oplus M'$ by requiring that $(T \perp T')(x+x')=T(x)+T(x')$ and $(Q \perp Q')(x+x')=Q(x)+Q(x')$.  (Note that $T(x,x)=2Q(x)$ for all $x \in M$ so if $2 \neq 0 \in R$ then the second condition follows from the first.)

Let $Q:M \to R$ be a quadratic form and suppose that $M$ is free (of finite rank).  In this case, a basis $e_1,\dots,e_n$ for $M$ gives an isomorphism $M \cong R^n$ in which $Q$ can be written
\[ Q(x)=Q(x_1e_1+\dots+x_ne_n)=\sum_i Q(e_i)x_i^2 + \sum_{i<j} T(e_i,e_j)x_ix_j \]
with $x=(x_1,\dots,x_n) \in R^n$.

For $a \in R$, the quadratic form $Q(x)=ax^2$ on $R$ is denoted $\la a \ra$; similarly, for $a_1,\dots,a_n \in R$, we abbreviate $\la a_1 \ra \perp \dots \perp \la a_n \ra = \la a_1,\dots,a_n \ra$.  For $a,b,c \in R$, the quadratic form $Q(x,y)=ax^2+bxy+cy^2$ on $R^2$ is denoted $[a,b,c]$.  

\begin{exm} \label{quatisnonsing}
Let $B=\quat{a,b}{F}$ be a quaternion algebra over $F$.  Then as in Example \ref{quatstdinvo}, in the basis $1,i,j,ij$ we have $\nrd \cong \la 1,-a,-b,ab \ra \cong \la 1,-a \ra \perp -b \la 1,-a \ra$ if $\opchar F \neq 2$ and $\nrd \cong [1,1,a] \perp b[1,1,a]$ if $\opchar F = 2$.  

Similarly, for $\nrd_0:B_0 \to F$ we have $\nrd_0 \cong \la -a,-b,ab \ra \cong \la -a \ra \perp -b \la 1,-a \ra$ if $\opchar F \neq 2$ and $\nrd_0 \cong \la 1 \ra \perp b[1,1,a]$ if $\opchar F = 2$.
\end{exm}

\subsection*{Quadratic forms over DVRs}

Now let $R$ be a local PID.  Then $R$ has valuation $\ord_v:R \to \Z_{\geq 0} \cup \{\infty\}$ and uniformizer $\pi$.  If $R=F$ is a field, then $\pi=1$ and the valuation is trivial, i.e.\ $\ord_v(x)=0$ for $x \in F^\times$ (and $\ord_v(0)=\infty$).

Let $Q:M \to R$ be a quadratic form over $R$.  Then since $R$ is a PID, $M$ is free; let $n$ be the rank of $M$ over $R$.  We will now seek to find a basis for $R^n$ in which a quadratic form $Q$ has a particularly simple form: we will seek to diagonalize $Q$ as far as possible.  In cases where $2 \in R^\times$, we can accomplish a full diagonalization; otherwise, we can at least break up the form orthogonally into indecomposable and distinguished forms of dimension at most $2$, as follows.

A quadratic form $Q$ over $R$ is \emph{atomic} if either:
\begin{enumroman}
\item $Q \cong \la a \ra$ for some $a \in R^\times$, or
\item $2 \not \in R^\times$ and $Q \cong [a,b,c]$ with $a,b,c \in R$ satisfying 
\[ \ord_v(b) < \ord_v(2a) \leq \ord_v(2c)\text{ and }\ord_v(a)\ord_v(b)=0. \] 
\end{enumroman}
In case (ii), we necessarily have $\ord_v(2)>0$ and $\ord_v(b^2-4ac)=2\ord_v(b)$.  

\begin{exm}
If $2 \in R^\times$, then a quadratic form $Q$ is atomic if and only if $Q(x)=ax^2$ for $a \in R^\times$.  
\end{exm}

\begin{exm} \label{opchar2atomic}
If $R=F$ is a field with $\opchar F=2$, then $[a,b,c]$ is atomic if and only if $b \in F^\times$; scaling $y$ by $a/b$ realizes this form as isomorphic to $a[1,1,ca/b^2]$ with $a \in F^\times$.  Therefore, over fields, recording the middle coefficient is unnecesary, and indeed other texts use $[a,b]$ to denote the quadratic form $ax^2+xy+by^2$.  

For example, take $R=\Z_2[\sqrt{2}]$ with normalized valuation $\ord_v(\sqrt{2})=1$ and let $Q(x,y)=x^2+\sqrt{2}xy$.  Then according to our definition, $Q$ is atomic, since $\ord_v(b)=1<\ord_v(2a)=2 \leq \ord_v(2c)=\infty$ and $\ord_v(a)=0$.  But this form is not globally divisible by any element of positive valuation, and a calculation shows that any isomorphic (equivalent) form has middle coefficient of positive valuation.
\end{exm}

\begin{exm}
Suppose $R=\Z_2$ is the ring of $2$-adic integers, so that $\ord_v(x)=\ord_2(x)$ is the largest power of $2$ dividing $x \in \Z_2$.  Recall that $\Z_2^\times/\Z_2^{\times 2}$ is represented by the elements $\pm 1, \pm 5$, therefore a quadratic form $Q$ over $\Z_2$ is atomic of type (i) above if and only if $Q(x) \cong \pm x^2$ or $Q(x) \cong \pm 5x^2$.  For forms of type (ii), the conditions $\ord_v(b) < \ord_v(2a) = \ord_v(a)+1$ and $\ord_v(a)\ord_v(b)=0$ imply in fact $\ord_v(b)=0$, and so a quadratic form $Q$ over $\Z_2$ is atomic of type (ii) if and only if $Q(x,y) \cong ax^2 + xy + cy^2$ with $\ord_2(a) \leq \ord_2(c)$.  Replacing $x$ by $ux$ and $y$ by $u^{-1} y$ for $u \in \Z_2^\times$ we may assume $a$ is a power of $2$, and then the atomic representative $[2^t, 1, c]$ of the isomorphism class of $Q$ is unique.
\end{exm}

A quadratic form $Q$ is \emph{decomposable} if $Q$ can be written as the orthogonal sum of two quadratic forms ($Q \cong Q_1 \perp Q_2$) and is \emph{indecomposable} otherwise.  

It follows by induction on the rank of $M$ that $Q$ is the orthogonal sum of indecomposable forms.  We will soon give an algorithmic proof of this fact and write each indecomposable form as a scalar multiple of an atomic form.  We begin with the following lemma.

\begin{lem}
An atomic form $Q$ is indecomposable.
\end{lem}

\begin{proof}
If $Q$ is atomic of type (i) then the space underlying $Q$ has rank $1$, so this is clear.  So suppose $Q=[a,b,c]$ is atomic of type (ii) and suppose $Q$ is decomposable.  It follows that if $x,y \in M$ then $T(x,y) \in 2R$.  Thus we cannot have $\ord_v(b)=0$, so $\ord_v(a)=0$, and further $\ord_v(b) \geq \ord_v(2)=\ord_v(2a)$; this contradicts the fact that $Q$ is atomic.
\end{proof}

\begin{prop} \label{qformdiag}
Let $R$ be a local PID and let $Q:M \to R$ be a quadratic form.  Then there exists a basis of $M$ such that the form $Q$ can be written
\[ Q \cong \pi^{e_1}Q_1 \perp \dots \perp \pi^{e_n} Q_n \]
where the forms $Q_i$ are atomic and $0 \leq e_1 \leq \dots \leq e_n \leq \infty$.
\end{prop}

In the above proposition, we interpret $\pi^{\infty} = 0$.  A form as presented in Proposition \ref{qformdiag} is called \emph{normalized}, and the integer $e_i$ is called the \emph{valuation} of $\pi^{e_i} Q_i$.  The tuple of valuations $e_i$ for $Q$ is unique.

\begin{exm}
By Example \ref{quatisnonsing}, if $B$ is a quaternion algebra over a field $F$ then the quadratic form $\nrd$ is normalized in the basis $1,i,j,ij$, with a similar statement for $\nrd_0$.
\end{exm}

We give an algorithmic proof of Proposition \ref{qformdiag}.  (Over fields, see Lam \cite[\S 1.2]{Lam}, and see Scharlau \cite[\S 9.4]{Scharlau} for fields of characteristic $2$.)

\begin{alg} \label{normalizequadform}
Let $R$ be a computable ring which is a local PID with (computable) valuation $\ord_v:R \to \Z_{\geq 0} \cup \{\infty\}$.

Let $Q:M \to R$ be a quadratic form over $R$ and let $e_1,\dots,e_n$ be a basis for $M$.  This algorithm returns a basis of $M$ in which $Q$ is normalized.
\begin{enumalg}
\item If $T(e_i,e_j)=0$ for all $i,j$, return $f_i := e_i$.  Otherwise, let $(i,j)$ with $1 \leq i \leq j \leq n$ be such that $\ord_v T(e_i,e_j)$ is minimal, taking $i=j$ if possible and if not taking $i$ minimal.
\item If $i=j$, let $f_1 := e_i$ and proceed to Step 3.  If $i \neq j$ and $2 \in R^\times$, let $f_1 := e_i + e_j$ and proceed to Step 3.  Otherwise, proceed to Step 4.
\item Let $e_i := e_1$.  For $k=2,\dots,n$ let 
\[ f_k := e_k-\frac{T(f_1,e_k)}{T(f_1,f_1)} f_1. \]
Let $m=2$ and proceed to Step 5.
\item (We have $2 \not\in R^\times$ and $i \neq j$.)  Let 
\[ f_1 := \frac{\pi^{\ord_v T(e_i,e_j)}}{T(e_i,e_j)} e_i, \] 
$f_2 := e_j$, $e_i := e_1$ and $e_j := e_2$.  Let $d := T(f_1,f_1)T(f_2,f_2)-T(f_1,f_2)^2$.  For $k=3,\dots,n$, let
\begin{align*}
t_k &:= T(f_1,f_2)T(f_2,e_k)-T(f_2,f_2)T(f_1,e_k) \\
u_k &:= T(f_1,f_2)T(f_1,e_k)-T(f_1,f_1)T(f_2,e_k)
\end{align*}
and let
\[ f_k := e_k+\frac{t_k}{d}f_1+\frac{u_k}{d} f_2. \]
Let $m=3$.
\item Recursively call the algorithm with $M=Rf_m \oplus \dots \oplus Rf_n$, and return $f_1,\dots,f_{m-1}$ concatenated with the output basis.
\end{enumalg}
\end{alg}

Given such a basis, one recovers the normalized quadratic form by factoring out in each atomic form the minimal valuation achieved.  (One can also keep track of this valuation along the way in the above algorithm, if desired.)

\begin{rmk} \label{2inrtimes}
Note that if $2 \in R^\times$, then this algorithm computes a diagonalization of the form $Q$, ordering the coefficients by their valuation.
\end{rmk}

\begin{proof}[Proof of correctness]
In Step 3, we verify that $\ord_v T(f_1,f_1) \leq \ord_v T(f_1, e_k)$.  Indeed, we have  
\[ T(f_1,f_1)=T(e_i,e_i)+2T(e_i,e_j) + T(e_j,e_j) \]
and so $\ord_v T(f_1,f_1)=\ord_v T(e_i,e_j)$ by the ultrametric inequality and the hypotheses that $\ord_v T(e_i,e_j) < \ord_v T(e_i,e_i),T(e_j,e_j)$ and $\ord_v(2)=0$.  So Steps 2 and 3 give correct output.

We have left to check Step 4.  This is proven by letting $f_k = e_k + t_k f_1 + u_k f_2$ and solving the linear equations $T(f_1,f_k) = T(f_2,f_k)=0$ for $t_k,u_k$.  The result then follows from a direct calculation, coupled with the fact that $\ord_v(d) = 2\ord_v T(f_1,f_2) \leq \ord_v(t_k)$ (and similarly with $u_k$).  This case only arises if (and only if)
\[ \ord_v T(f_1,f_2) < \ord_v T(f_1,f_1) = \ord_v(2Q(f_1)) \leq \ord_v(2Q(f_2)) \]
so the corresponding block is indeed atomic.
\end{proof}

\begin{exm}
Consider the binary quadratic form $[a,b,c]$ over $\Z_2$.  Then $T(e_1,e_1)=2a$, $T(e_1,e_2)=b$, and $T(e_2,e_2)=2c$.  We follow the course of Algorithm \ref{normalizequadform}.  If $\ord_v(2a)$ is minimal, then in Steps 2 and 3 we diagonalize (complete the square): we have $f_1=e_1$ and $f_2=e_2-(b/2a)e_1$ and so we obtain the (isomorphic) form $\la a, c+b^2/4a \ra$.  If $\ord_v(2c)$ is minimal, then we similarly obtain $\la c, a+b^2/4c \ra$.  Finally, if $\ord_2(b)$ is minimal, then we enter Step 4.  Since $(i,j)$ was taken with $i$ minimal, for illustration we may suppose $i=1$ and $j=2$.  Then we have $t=\ord_v(b)<\ord_v(2a) \leq \ord_v(2c)$.  Writing $a=2^t a'$, $b'=2^t b'$ and $c'=2^t c'$, in Step 4, we simply have $f_1=(1/b')e_1$ and $f_2=e_2$ and we obtain the form $2^t [a'/(b')^2, 1, c']$ and $[a'/(b')^2,1,c']$ is indeed atomic.  
\end{exm}

\begin{exm}
Consider the form $q(x,y,z)=xy+xz$ over $\Z_2$.  We enter Step 4 with $f_1 = e_1$ and $f_2=e_2$.  We compute that $d=-T(f_1,f_2)=-1$, and $t_3=0$ and $u_3=1$.  Thus $f_3=e_3-f_2=e_3-e_2$, and we obtain the form $[0,1,0] \perp \la 0 \ra$.
\end{exm}

We note that Algorithm \ref{normalizequadform} requires $O(n^2)$ arithmetic operations in $R$.  This algorithm can be modified suitably to operate on the Gram matrix $(T(e_i,e_j))_{i,j}$ of the quadratic form $Q$, which as explained above recovers the quadratic form when $2 \neq 0 \in R$.  

For a quadratic form $Q:M \to R$, we define
\[ \rad(Q)=\{x \in M: T(x,y)=0 \text{ for all $y \in M$}\}; \]
we say $Q$ is \emph{nonsingular} if $\rad(Q)=\{0\}$.  

\begin{exm} \label{exmisnonsing}
We have $\rad(Q \perp Q')=\rad(Q) \oplus \rad(Q')$, and if $Q$ is atomic then $\rad(Q)=\{0\}$.  In particular, one can read off $\rad(Q)$ directly from a normalized form by the corresponding valuations.
\end{exm}

\subsection*{Identifying quaternion algebras}

Using the above normalization of a quadratic form in the case where $R=F$ is a field, we can directly identify quaternion algebras amongst algebras with a standard involution.

\begin{prop} \label{identquat}
Let $B$ be an $F$-algebra with a standard involution.  If $\dim_F B=4$, then $B$ is a quaternion algebra if and only if $\nrd$ is nonsingular.
\end{prop}

\begin{proof}
If $B$ is a quaternion algebra, then $\nrd$ is nonsingular by Example \ref{quatisnonsing}.

Conversely, $B$ has a basis $1,i,j,k$ which is a normalized basis for $Q$.  First suppose $\opchar F \neq 2$.  By orthogonality we have $\trd(i)=0$ so $i^2=-\nrd(i)=a \neq 0$ by nonsingularity and similarly $j^2=b \neq 0$, and $ji+ij=0$ from (\ref{Bfm}) so $(ij)^2=-ab$.  Thus $B \supset \quat{a,b}{F}$ hence this map is an isomorphism.  The case $\opchar F = 2$ follows similarly: now instead we have $i^2+i=a$ and $ji=\overline{i}j=(i+1)j$.  
\end{proof}

Proposition \ref{identquat} yields the following algorithm.

\begin{alg} \label{identquatalg}
Let $B$ be an $F$-algebra with $\dim_F B=4$ (specified by a multiplication table).  This algorithm returns \textsf{true} if and only if $B$ is a quaternion algebra, and if so returns an isomorphism $B \cong \quat{a,b}{F}$.

\begin{enumalg}
\item Verify that $B$ has a standard involution by calling Algorithm \ref{deg2}.  If not, return \textsf{false}.
\item Compute a normalized basis $1,i,j,k$ for the quadratic form $\nrd:B \to F$ by calling Algorithm \ref{normalizequadform}.  
\item Test if $\nrd$ is nonsingular as in Example \ref{exmisnonsing}.  If so, return \textsf{true} and the quaternion algebra $\quat{a,b}{F}$ given by the standard generators $i,j$.
\end{enumalg}
\end{alg}

\begin{rmk}
Given a quaternion algebra over $\Q$, R\'onyai \cite[Theorem 2.1]{Ronyai1} gives an algorithm to compute a standard representation, but this algorithm tests a polynomial of degree $2$ over $\Q$ for irreducibility; the above algorithm requires no such test.
\end{rmk}

\begin{rmk}
If in Step 3 one finds that $\nrd$ is not nonsingular, then one has the further refinement of Algorithm \ref{identquatalg} as follows.

We denote by $\rad(B)$ the \emph{Jacobson radical} of $B$, the largest two-sided \emph{nil ideal} of $B$, i.e.\ the largest two-sided ideal in which every element is nilpotent.  An algebra $B$ for which $\rad(B)=\{0\}$ is called \emph{semisimple}.
We claim that $\rad(B)=\rad(\nrd)$.  Indeed, let $e \in B$ be nilpotent, so that $e^2=0$.  For any $x \in B$, we have by (\ref{Bfm}) that
\[ xe+ex=\trd(x)e+\trd(xe). \]
It follows that $e$ generates a nil ideal if and only if $T(x,e)=0$ for all $x \in B$, which holds if and only if $x \in \rad(\nrd)$.  Thus $\rad(B)=\rad(\nrd)$.  One can then easily modify the algorithm to output $\rad(B)=\rad(\nrd)$.  
\end{rmk}

\begin{rmk}
Another algorithm which tests if $B$ is a quaternion algebra (but does not give a standard representation) under the assumption $\opchar F = 0$ runs as follows.  (See Lam \cite[Chapter 4]{Lam} for the standard facts we use.)  By the Wedderburn-Artin theorem and a dimension count, the algebra $B$ over $F$ is a quaternion algebra if and only if $B$ is central and semisimple.  We verify that $B$ is central as in Remark \ref{computecenter}.  To verify semisimplicity, if $\opchar F=0$, Dickson \cite[\S 66]{Dickson} showed that $B$ with $\dim_F B=n$ is semisimple if and only if the matrix $(\Tr(e_ie_j))_{i,j=1,\dots,n}$ has full rank $n$, where $\Tr$ is the (left) algebra trace.
\end{rmk}

In view of Algorithm \ref{identquatalg}, we assume from now on that a quaternion algebra $B$ over a field $F$ is given as input by a standard representation.

Over a general domain $R$, the above algorithms do not generalize directly, as we cannot hope to normalize a quadratic form in such a simple way for over rings that are no longer local PIDs.  Indeed, the category of quadratic forms over a general domain $R$ can be quite complicated---already forms over the integers $\Z$ are of significant interest.  However, over Dedekind domains, we can still recognize quaternion orders, and one instead understands these orders as in Section 1 via their localizations, a subject which will consume the later sections of this article.

\subsection*{Identifying quaternion orders}

Let $F$ be a number field and let $\Z_F$ be its ring of integers.  In this section, we give an algorithm which allows us in many cases to put quaternion orders in a standard form as in the discussion of Lemma \ref{SibR}.

\begin{alg} \label{identquatorder}
Let $\calO \subset B$ be a quaternion order over $\Z_F$.  Let $\iota:K \to B$ be an embedding of $F$-algebras with $K$ a field such that $[K:F]=2$ and $\iota(K) \cap \calO = \Z_K$ is maximal.  This algorithm returns a fractional ideal $\frakb$ of $K$, an element $j \in \calO$ such that $B = \iota(\Z_K) \oplus \iota(\frakb) j \cong \quat{\Z_K,\frakb,b}{\Z_F}$.

\begin{enumalg}
\item Identify $K$ with $\iota(K)$.  Let $K=F \oplus Fi$ with $i \in B$.  Compute $j \in B$ orthogonal to $1,i$.  
\item Let $x_1,\dots,x_m$ be a generating set for $\calO$ as a $\Z_F$-module.  Write $x_k=a_k + b_k j$ with $a_k,b_k \in K$ for $k=1,\dots,m$. 
\item Compute a pseudo-basis $\Z_K \oplus \frakb j$ for the $\Z_K$-module generated by $(a_k,b_k)$ for $k=1,\dots,m$ using a HMF.
\item Let $a,b$ be generators for $\frakb$ as an $\Z_F$-module.  If $\trd(j) \neq 0$, then let $c := \trd(bj)a - \trd(aj)b$, let $j := cj$ and $\frakb := (1/c)\frakb$.  Return $\frakb$ and the element $j$.
\end{enumalg}
\end{alg}

\begin{proof}[Proof of correctness]
In Step 4, we check directly that $\trd(j)=\trd(ij)=0$, as desired.
\end{proof}

\begin{rmk}
One can extend Algorithm \ref{identquatorder} when $\iota(K) \cap \calO=S$ is no longer maximal by an appropriate modification of the HMF algorithm over $S$.
\end{rmk}

\section{Identifying the matrix ring}

In this section, we continue the pursuit of our motivating question and address the computational complexity of identifying the matrix ring over a field.  Throughout this section, let $F$ be a computable field.  We represent a quaternion algebra $B$ over $F$ by a standard form $B=\quat{a,b}{F}$.

\begin{prob*}[\textsf{IsMatrixRing}]
Given a quaternion algebra $B$ over $F$, determine if $B \cong \M_2(F)$.
\end{prob*}

We may also ask for a solution to the more difficult problem of constructing an explicit isomorphism.

\begin{prob*}[\textsf{ExhibitMatrixRing}]
Given a quaternion algebra $B$ over $F$, determine if $B \cong \M_2(F)$ and, if so, output such an isomorphism.
\end{prob*}


\subsection*{Zerodivisors}

Let $B$ be a quaternion algebra.  The following structural lemma allows us to address the above problems.

\begin{lem} \label{identm2f}
The following are equivalent:
\begin{enumroman}
\item $B \cong \M_2(F)$;
\item $B$ is not a division ring;
\item There exists a nonzero $e \in B$ such that $e^2=0$; and
\item $B$ has a proper, nonzero left (or right) ideal $I$.
\end{enumroman}
\end{lem}

If $B \cong \M_2(F)$, we say that $B$ is \emph{split}.  More generally, if $K \supset F$ is a field containing $F$, then we say $K$ is a \emph{splitting field} for $B$ if $B_K=B \otimes_F K$ is split.

We give a proof of Lemma \ref{identm2f} in an algorithmically effective way in this section.  The implication (i) $\Rightarrow$ (ii) is clear.  The implication (ii) $\Rightarrow$ (iii) is obtained as follows.

\begin{alg} \label{etatoeps}
Let $x \in B$ be a zerodivisor.  This algorithm returns a nonzero element $e \in B$ such that $e^2=0$.

\begin{enumalg}
\item If $\trd(x)=0$, return $x$.
\item Compute $0 \neq y \in B$ orthogonal to $1,x$ with respect to the quadratic form $\nrd$.  If $xy=0$, return $y$; otherwise, return $xy$.
\end{enumalg}
\end{alg}

\begin{proof}[Proof of correctness]
The element $x \neq 0$ is a zerodivisor if and only if $\nrd(x)=x\overline{x}=0$.  Since $y$ is orthogonal to $1$ we have $\trd(y)=0$ so $\overline{y}=-y$; similarly, since $y$ is orthogonal to $x$ we have $\trd(xy)=-\trd(x\overline{y})=0$.  If $xy=0$ then $y$ is a zerodivisor.  If $xy \neq 0$ then $\nrd(xy)=\nrd(x)\nrd(y)=0$, as desired.  
\end{proof}

The implication (iii) $\Rightarrow$ (iv) follows, since $e$ generates a proper left (or right) ideal.  Below, in the proof of correctness of the following algorithm, we will show that if $I=Be$ then $\dim_F I=2$; the final implication (iv) $\Rightarrow$ (i) then follows since left multiplication gives a nonzero $F$-algebra map $B \to \End_F(I) \cong \M_2(F)$ which is injective since $B$ is simple and therefore an isomorphism as $\dim_F B = 4 = \dim_F \M_2(F)$.  

\begin{alg} \label{m2ffound}
Let $e \in B$ satisfy $e^2=0$.  This algorithm returns a standard representation $B \cong \quat{1,1}{F} \cong \M_2(F)$.

\begin{enumalg}
\item Find $k \in \{i,j,ij\}$ such that $\trd(ek)=s \neq 0$.  
Let $t=\trd(k)$ and $n=\nrd(k)$, and let $e'=(1/s)e$.
\item Let $j'=k+(-tk+n+1)e'$ and let
\[ i' = \begin{cases}
e'k-(k+t)e', & \text{ if $\opchar F \neq 2$}; \\
k+((t+1)k+n+1)e', & \text{ if $\opchar F = 2$}.
\end{cases} \]
Return $i',j'$.
\end{enumalg}
\end{alg}

\begin{proof}[Proof of correctness]
In Step 1, if $\trd(ek)=0$ for all such $k$ then $e \in \rad(\nrd)$, contradicting Lemma \ref{identquat}.  We have $\trd(e'k)=\trd(ke')=1$ so $\trd(\overline{e'}k)=-1$.

Consider $I=Fe' + Fke'$.  Note $\trd(ke') \neq 0$ implies that $e',ke'$ are linearly independent.  Let $A$ be the subalgebra of $B$ generated by $e'$ and $k$.  We have $e'k+ke'=te'+1$ from (\ref{Bfm}) and $k^2=tk-n$, and thus we compute that left multiplication yields a map
\begin{align*}
A &\to \End_F(I) \cong \M_2(F) \\ 
e',k &\mapsto \begin{pmatrix} 0 & 1 \\ 0 & 0 \end{pmatrix}, \begin{pmatrix} 0 & -n \\ 1 & t \end{pmatrix}.
\end{align*}
A direct calculation then reveals that $j' \mapsto \begin{pmatrix} 0 & 1 \\ 1 & 0 \end{pmatrix}$ and $i' \mapsto \begin{pmatrix} 1 & 0 \\ 0 & -1 \end{pmatrix}$ if $\opchar F \neq 2$ and $i' \mapsto \begin{pmatrix} 0 & 1 \\ 1 & 1 \end{pmatrix}$ if $\opchar F = 2$, as in Example \ref{M2Fisom11}.

It follows all at once that $A=B$, that $I=Be'$, and that the map $B \to \M_2(F)$ is an isomorphism.
\end{proof}

\begin{rmk}
An algorithm like the above which requires linear algebra in $F$ is claimed but not exhibited explicitly by R\'onyai \cite{Ronyai1}; see also further of R\'onyai \cite[\S 5.1]{RonyaiExplicitIsomorphism}.
\end{rmk}

\subsection*{Conics} \label{conics}

We have already seen in Lemma \ref{identm2f} that $B \cong \M_2(F)$ if and only if there exists $0 \neq e \in B$ such that $e^2=0$.  To this end, as in the previous section let 
\[ B_0=\{e \in B:\trd(e)=0\}. \]
We have $\dim_F B_0=3$, and given a standard representation for $B=\quat{a,b}{F}$, we have a basis for $B_0$ given by $i,j,ij$ if $\opchar F \neq 2$ and $1,j,ij$ if $\opchar F=2$, as in Example \ref{quatisnonsing}.  

We may identify the set $\PP(B_0)=B_0^\times/F^\times$ with the points of the projective plane $\PP^2(F)$ over $F$.  Then the equation $\nrd_0(x,y,z)=0$ yields a \emph{conic} $C \subset \PP_F^2$ defined over $F$, a nonsingular projective plane curve of degree $2$.

\begin{lem} \label{identm2fconic}
The following are equivalent:
\begin{enumroman}
\item $B \cong \M_2(F)$;
\item[(v)] The quadratic form $Q=\nrd|_{B_0}$ associated to $B$ represents zero over $F$; and
\item[(vi)] The conic $C$ associated to $B$ has an $F$-rational point.
\end{enumroman}
\end{lem}

Therefore we are led to the following problems.

\begin{prob}[\textsf{HasPoint}]
Given a conic $C$ defined over a field $F$, determine if $C$ has an $F$-rational point.
\end{prob}

\begin{prob}[\textsf{ExhibitPoint}]
Given a conic $C$ defined over a field $F$, determine if $C$ has an $F$-rational point and, if so, output such a point.
\end{prob}

These problems could be equivalently formulated as follows: given a nonsingular ternary quadratic form $Q:V \to F$, determine if $F$ is \emph{isotropic} (represents zero nontrivially) and, if so, find $0 \neq x \in V$ such that $Q(x)=0$.  We find the geometric language here to be more suggestive, but really these are equivalent ways to describe the same situation.

By Algorithm \ref{normalizequadform}, given a conic $C$ over $F$, there is a (deterministic, polynomial-time) algorithm which computes a change of coordinates in which $C$ is given by the equation
\[ ax^2+by^2+cz^2=0 \]
if $\opchar F \neq 2$, with $a,b,c \in F^\times$, and
\[ ax^2+axy+aby^2+cz^2=0 \]
if $\opchar F=2$, with $a,c \in F^\times$ and $b \in F$ by Example \ref{opchar2atomic}.  In the first case, multiplying through by $abc \neq 0$ we obtain $bc(ax)^2+ac(by)^2+(abc^2)z^2=0$ which arises as the form associated to $\quat{-bc,-ac}{F}$; in the second case, we multiply through by $c \neq 0$ to obtain $(ac)x^2+(ac)xy+b(ac)y^2+(cz)^2=0$ which is associated to $\quat{b,ac}{F}$.  
Together with Algorithm \ref{m2ffound}, therefore, we arrive at the following lemma.

\begin{prop}
The association $B \mapsto C=\nrd_0$ gives a bijection between quaternion algebras over $F$ up to isomorphism and conics over $F$ up to isomorphism.  

Problems \probsf{IsMatrixRing}, \probsf{ExhibitMatrixRing} are (deterministic poly\-nomial-time) equivalent to Problems \probsf{HasPoint}, \probsf{ExhibitPoint}, respectively.
\end{prop}

\begin{proof}
We need only identify isomorphisms: we need to show that two quaternion algebras $B \cong B'$ are isomorphic if and only if the induced conics $C \cong C'$ are isomorphic. 

We treat only the case $\opchar F \neq 2$; the case $\opchar F = 2$ follows similarly.   If $\phi:B \to B'$ is an isomorphism of quaternion algebras, then $\phi(1)=1$ so $\phi(B_0)=B'_0$, and the reduced norm is determined by the standard involution which is unique, so $\nrd_B = \nrd_{B'} \circ \phi$.  

Conversely, suppose $\psi:C \to C'$ is an isomorphism.  Choose a quadratic form $Q$ so that $C$ is given by $Q=0$ in $\PP_F^2$, normalized and scaled so that $Q \cong \nrd_0$ for some $B \cong \quat{a,b}{F}$.  Choose similarly $Q'$ for $C'$.  Then $\psi$ is given by an element of $\PGL_3(F)$ and there exists a lift of $\psi$ to $\GL_3(F)$ such that $Q = Q' \circ \psi$.  The $F$-linear map $\psi:B_0 \to B'_0$ extends naturally (defining $\phi(1)=1$) to an $F$-linear map which we also denote $\psi:B \to B'$, and we must show that $\psi$ is an $F$-algebra isomorphism.  

Suppose $B=\quat{a,b}{F}$.  Then we have $\nrd(\psi(i))=\nrd(i)=-a$ and $\nrd(\psi(i))=\psi(i)\overline{\psi(i)}=-\psi(i)^2$ so $\psi(i)^2=a$.  Similarly we have $\psi(j)^2=b$.  We have $ji=-ij$ since $i,j$ are orthogonal, but then $\psi(i),\psi(j)$ are orthogonal so $\psi(j)\psi(i)=-\psi(i)\psi(j)$.  Finally, we have that both $\psi(ij)$ and $\psi(i)\psi(j)$ are orthogonal to $1,\psi(i),\psi(j)$, and $\psi(ij)^2=-ab=(\psi(i)\psi(j))^2$, so $\psi(ij)=\pm \psi(i)\psi(j)$.  If the negative sign occurs, we replace $\psi$ by the linear map defined on the basis $1,i,j,ij$ unmodified on $1,i,j$ but negated on $ij$; this map is now an $F$-algebra homomorphism.  Together, these imply that $B' \cong \quat{a,b}{F}$ as well.  
\end{proof}

We conclude this section by considering a simple case of the above problems.  First, let $F=\F_q$ be a finite field with $q$ elements.  Indeed, Problem \probsf{HasPoint} is trivial: since every conic over a finite field has a point (an elementary argument), one can simply always output \textsf{true}!

For problem (\textsf{ExhibitPoint}), we will make use of the following related problem.

\begin{prob}[\textsf{SquareRoot}]
Given $a \in F^{\times 2}$, output $b \in F^\times$ such that $b^2=a$.
\end{prob}

We have two cases.  First, if $q$ is even, then one can solve Problem \probsf{SquareRoot} in deterministic polynomial time (by repeated squaring, since $q-1=\#\F_{2^r}^\times$ is odd); for a conic in the form given in Example \ref{quatisnonsing}, given up to scaling by $x^2+by^2+byz+abz^2$ with $a,b \in \F_q$ and $b \neq 0$, this is already sufficient to solve Problem \probsf{ExhibitPoint}.  If $q$ is odd, then there exists a deterministic polynomial-time algorithm to solve (\textsf{ExhibitPoint}) over $\F_q$ by work of van de Woestijne \cite{Christiaan}.  There also exists a probabilistic polynomial-time algorithm, which intersects the conic with a random line and then calls (\textsf{SquareRoot}), and there is a probabilistic polynomial-time algorithm to solve (\textsf{SquareRoot}) but no deterministic such algorithm (without further assumption of a generalized Riemann hypothesis).  The latter algorithm is extremely efficient in practice.  

\begin{rmk}
It would also be interesting to study the corresponding problem where $\M_2(F)$ is replaced by another quaternion algebra $B'$: in other words, to test if two quaternion algebras $B$, $B'$ over $F$ are isomorphic and, if so, to compute an explicit isomorphism.  Since the reduced norm is determined by the standard involution on a quaternion algebra, and this involution is unique, it follows that if $B \cong B'$ then $\nrd_B \cong \nrd_{B'}$; in fact, this is an equivalence even when restricted to the trace zero subspace \cite{Lam}.  Therefore one is led to consider the problem of determining if two quadratic forms are isometric and, if so, to compute an explicit isometry.  
\end{rmk}

\begin{rmk}
More generally, one can establish a functorial bijection between twisted similarity classes of ternary quadratic forms over a commutative ring $R$ and quaternion rings over $R$ via the Clifford algebra; see work of the author \cite{VoightCrelle}.  It would be interesting to investigate the algorithmic implications of this correspondence.
\end{rmk}

\section{Splitting fields and the Hilbert symbol}

In this section, we exhibit algorithms for solving the Problem \probsf{IsMatrixRing} over a local field with residue characteristic not 2: in this setting, our problem is otherwise known as computing the Hilbert symbol.

\subsection*{Hilbert symbol}

Let $F$ be a field with $\opchar F \neq 2$, and let $a,b \in F^\times$.  The \emph{Hilbert symbol} is defined to be
\[ (a,b)_F = 
\begin{cases}
1, &\text{ if $\quat{a,b}{F} \cong \M_2(F)$;} \\
-1, &\text{ otherwise.}
\end{cases} \]

We begin by recalling a well-known criterion \cite[Corollaire 2.4]{Vigneras}.

\begin{lem} \label{isnormdet}
A quaternion algebra $\quat{a,b}{F}$ is split if and only if $b \in N_{K/F}(K^\times)$, where $K=F[i]$.
\end{lem}

Here, we write $K=F[i]=F \oplus Fi$ to be the quadratic $F$-algebra generated by $i$.

\begin{proof}
If $\N_{K/F}(u+vi)=\nrd(u+vi)=b$ with $x,y \in F$, then $x=u+vi+j$ has $\nrd(x)=\nrd(u+vi+j)=\nrd(u+vi)+\nrd(j)=b-b=0$, so $B$ is not a division ring, so $B \cong \M_2(F)$ by Lemma \ref{identm2f}.  Conversely, if $B \xrightarrow{\sim} \M_2(F)$, then after conjugating by an element of $\GL_2(F)$ we may assume $i \mapsto \begin{pmatrix} 0 & a \\ 1 & 0 \end{pmatrix}$ (rational canonical form).  The condition that $ji=-ij$ implies that $j \mapsto \begin{pmatrix} u & -av \\ v & -u \end{pmatrix}$ and $j^2=u^2-av^2=b=\N_{K/F}(u+vi)$.
\end{proof}

\begin{lem}
We have $(a,b)_F=(b,a)_F$ and $(a,b)_F=(-ab,b)_F$.  If $u,v \in F^\times$ then $(a,b)_F=(au^2,bv^2)_F$.
\end{lem}

\begin{proof}
Interchanging $i,j$ gives an isomorphism $\quat{a,b}{F} \cong \quat{b,a}{F}$; replacing $i,j$ by $ui,vj$ gives an isomorphism $\quat{a,b}{F} \cong \quat{u^2a, v^2b}{F}$.  By considering the algebra generated by $ij,j$ we see that $\quat{a,b}{F} \cong \quat{a,-ab}{F}$.
\end{proof}

\subsection*{Local Hilbert symbol}

For the rest of this section, let $F$ be a number field.  For a place $v$ of $F$, let $F_v$ denote the completion of $F$ at $v$ and let $R_v$ be its valuation ring.  Let $\pi_v$ be a uniformizer for $F_v$ and let $k_v$ be the residue field of $F_v$.

If $a,b \in F_v^\times$, we abbreviate $(a,b)_v = (a,b)_{F_v}$.  We now proceed to discuss the computability of $(a,b)_v$, and thereby Problem \probsf{IsMatrixRing} for local fields $F_v$ with $\opchar k_v \neq 2$.

\begin{rmk} \label{normmapsurj}
With Lemma \ref{isnormdet} in mind, we recall the following facts about local norms.  There is a unique unramified quadratic extension $K_v$ of $F_v$, obtained from the corresonding unique such extension of residue fields.  Then $\N_{K_v/F_v}(K_v^\times) = R_v^\times \times \pi_v^{2\Z}$ by Hensel's lemma, since the norm map in an extension of finite fields is surjective.  For further details, see Neukirch \cite[Corollary V.1.2]{Neukirch} or Fr\"ohlich  \cite[Proposition 7.3]{Frohlich}.
\end{rmk}

We begin by recalling the following fundamental result concerning division quaternion algebras over a local field {\cite[Th\'eor\`emes II.1.1, II.1.3]{Vigneras}}.

\begin{lem} \label{vignlem}
Let $v$ be a noncomplex place of $F$.  Then there is a unique quaternion algebra $B_v$ over $F_v$ which is a division ring, up to $F_v$-algebra isomorphism.
\end{lem}

Note that there is no division quaternion algebra over $\C$ since $\C$ is algebraically closed.  The unique division algebra over $\R$ is the classical ring of Hamiltonians $\HH=\quat{-1,-1}{\R}$.  If $v$ is nonarchimedean, then the unique division ring over $F_v$ is given by $B_v \cong \quat{K_v,\pi_v}{F_v}$, where $K_v$ is the (unique) unramified quadratic extension of $F_v$.  


Let $B$ be a quaternion algebra over $F$.  We say $B$ is \emph{unramified} (or \emph{split}) at $v$ if $B \otimes_F F_v \cong \M_2(F_v)$, i.e.~$F_v$ is a splitting field for $B$; otherwise (if $B_v$ is a division ring) we say $B$ is \emph{ramified} at $v$.

A place $v$ of $F$ is \emph{odd} if either $v$ is real or $v$ is nonarchimedean and $\#k_v$ is odd; $v$ is \emph{even} if $v$ is nonarchimedean and $\#k_v$ is even.  (A complex place is neither odd nor even.)  For an odd place $v$ and $a \in F_v^\times$, we define the \emph{square symbol} 
\[ \issq{a}{v}=
\begin{cases}
1, & \text{ if $a \in F_v^{\times 2}$}; \\
-1, & \text{ if $a \not\in F_v^{\times 2}$ and $\ord_v(a)$ is even}; \\
0, & \text{ if $a \not\in F_v^{\times 2}$ and $\ord_v(a)$ is odd}.
\end{cases} \]
Here we set the convention that $v$ is a real place then $\pi_v=-1$ is a uniformizer for $F_v \cong \R$ and that $a=(-1)^{\ord_v(a)}|a|$; in other words, $\issq{a}{v}=1$ or $0$ according as $a>0$ or $a<0$.  

Suppose $v$ is nonarchimedean.  If $\ord_v(a)=0$, then $\issq{a}{v}=\legen{a}{v}$ is the usual Legendre symbol (see (\ref{Legendredef}) below); in fact, $\issq{a}{v}=0$ if and only if $\ord_v(a)$ is odd.  Note that the square symbol is not multiplicative, for example $\issq{\pi_v^2}{v}=1 \neq 0 = \issq{\pi_v}{v}^2$; it is multiplicative when restricted to the the subgroup of elements with even valuation, however.  

Finally, we note that $\issq{a}{v}=-1$ if and only if $F_v(\sqrt{a})$ is an unramified field extension of $F_v$ and $\issq{a}{v}=0$ if and only if $F_v(\sqrt{a})$ is ramified; when $v$ is real, we follow the convention that $\C$ is considered to be ramified over $\R$.

\begin{prop} \label{padicdet}
Let $v$ be an odd place of $F$ and let $a,b \in F_v^\times$.  Then $(a,b)_{v}=1$ if and only if
\[ \issq{a}{v}=1\quad\text{or}\quad \issq{b}{v}=1\quad\text{or}\quad \issq{-ab}{v}=1 \quad\text{or}\quad \issq{a}{v}=\issq{b}{v}=\issq{-ab}{v}=-1. \]
\end{prop}

\begin{proof}
First, suppose $v$ is archimedean.  Then $(a,b)_{v}=1$ if and only if $v(a)>0$ or $v(b)>0$ if and only if $\issq{a}{v}=1$ or $\issq{b}{v}=1$.  So we suppose $v$ is nonarchimedean.  

Let $B_v=\quat{a,b}{F_v}$, and let $K_v=F_v[i]$, where we recall $i^2=a$.  Since $(a,b)_v=(b,a)_v=(a,-ab)_v$, the statement is symmetric in interchanging $a,b$ and replacing $b$ by $-ab$.  If one of $\issq{a}{v}=1$ or $\issq{b}{v}=1$ or $\issq{-ab}{v}=1$, then we may suppose $\issq{a}{v}=1$; consequently, $K_v$ is not a field, so $B_v$ is not a division ring and by Lemma \ref{identm2f} we have $(a,b)_v=1$.  We cannot have $\issq{a}{v}=\issq{b}{v}=\issq{-ab}{v}=0$.  Thus we have only to consider the case $\issq{a}{v}=-1$.  

If $\issq{b}{v}=-1$, then since $K_v$ is the unique unramified quadratic extension of $F_v$ and $\ord_v(b)$ is even, we have $b \in \N_{K_v/F_v}(K_v^\times)$ by Remark \ref{normmapsurj}, so by Lemma \ref{isnormdet} we have that $B_v$ is split so $(a,b)_v=1$.  Otherwise, $\issq{b}{v}=0$.  But now $F_v[i]=K_v$ is the unramified quadratic extension of $F_v$ so $b \not\in \N_{K_v/F_v}(K_v^\times)$ and thus $B_v$ is a division ring by Lemma \ref{isnormdet}, so $(a,b)_v=-1$.
\end{proof}

\begin{cor} \label{legen}
Let $a,b \in R_v \setminus \{0\}$ and suppose $a \in R_v^\times$.  Then $(a,b)_v=\displaystyle{\legen{a}{v}^{\ord_v b}}$.
\end{cor}

\subsection*{Representing local fields} \label{localfields}

When discussing computability for local fields, we immediately encounter the following issue: a local field $F_v$ is uncountable, so it is not computable.

One has at least two choices for overcoming this obstacle.  One possibility is to use \emph{exact local field arithmetic}, where one includes with the specification of an element its precision.  One then requires the output of algorithms to be a continuous function of the input and to be correct with whatever output precision is given.  This way of working with $\R$ (or $\C$) also goes by the name \emph{exact real} (or \emph{complex}) \emph{arithmetic}.  This model has several advantages.  In practice, for many applications this works extremely well: if more precision is required in the output, one simply gives more precision in the input.  Consequently this model is also very efficient.  Although this method does not realize a local field $F$ as a computable field, all of the algorithms we discuss in this article work well in this model for $F_v$.

A second method is simply to work in a computable subfield $F$ of the local field $F_v$.  Indeed, any subfield $F$ which is countably generated over its prime field is computable.  In this article, we will take this approach; it is more appropriate for the theoretical discussion below (even as it will be less efficient in practice).

With this discussion in mind, we represent a local field as follows.  First, let $F$ be a number field.  Let $v$ be a place of $F$.  If $v$ is archimedean, then it is specified by some ordering of the roots of $f$ in $\C$.  If $v$ is nonarchimedean, then $v$ is specified by a prime ideal in the ring of integers in $F$.  We can thereby compute a uniformizer $\pi_v \in F$ for the place $v$ by the Chinese remainder theorem.  

We then represent the local field as $F_v^{\textup{alg}}=\overline{F} \cap F_v$, an algebraic closure of $F$ in $F_v$.  Given a (monic) polynomial $g$ with coefficients in $F$, there exists a deterministic algorithm which returns the roots of $g$ in $F_v$ (as elements of $F_v^{\textup{alg}}$).  In the nonarchimedean case, Hensel's lemma provides the essential ingredient to show that one can (efficiently) compute with $F_v^{\textup{alg}}$.  With this choice, by computing in the subfield generated by any element $x \in F_v^{\textup{alg}}$ we can compute the discrete valuation $\ord_v : F \to \Z \cup \{\infty\}$ as well as the reduction map $R_v \to k_v$ modulo $\pi_v$.  When $v$ is real, we recall that $\ord_v(a)=0,1$ according as $a>0$ or $a<0$, and so the computability of $\ord_v$ follows from well-known algorithms for exact real root finding.  

The above discussion applies equally well to the case of global function fields; see Remark \ref{globalfunfields}.  For more on computably algebraically closed fields, we refer again to Stoltenberg-Hansen and Tucker \cite{SHT}.

\subsection*{Computing the local Hilbert symbol}

To conclude, we discuss the computability of the Hilbert symbol for odd places using Proposition \ref{padicdet}.  We use Proposition \ref{padicdet} and the correspondence above to relate Problem (\textsf{HasPoint}) to the problem of computing the square symbol.  

Suppose $F_v$ is archimedean.  The Hilbert symbol for $F_v \cong \C$ is trivial.  If $v$ is real, then $\issq{a}{v}=1,0$ according as $a>0$ or $a<0$, so by the correspondence above this solves \probsf{HasPoint} for these fields.  It follows that Problem \probsf{ExhibitPoint} is equivalent to Problem \probsf{SquareRoot}, and there is a deterministic algorithm to solve this problem in the computable subfield $F_v^{\textup{alg}}=\overline{F} \cap \R$ by hypothesis.

Next, suppose $F_v$ is nonarchimedean and that $v$ is odd.  Then we can evaluate $\issq{a}{v}$ by simply computing $\ord_v(a)=e$; if $e$ is odd then $\issq{a}{v}=0$, whereas if $e$ is even then $\issq{a}{v}=\displaystyle{\legen{a_0}{v}}$ where $a_0=a \pi_v^{-e} \in R_v$ and $\displaystyle{\legen{a_0}{v}=\legen{a_0}{\frakp}}$ is the usual Legendre symbol, defined by
\begin{equation} \label{Legendredef}
\legen{a_0}{\frakp}=
\begin{cases}
0, &\text{ if $a_0 \equiv 0 \pmod{\frakp}$}; \\
1, &\text{ if $a_0 \not\equiv 0 \pmod{\frakp}$ and $a_0$ is a square modulo $\frakp$}; \\
-1, &\text{ otherwise}.
\end{cases}.
\end{equation}
The Legendre symbol can be computed in deterministic polynomial time by Euler's formula
\[ \legen{a_0}{\frakp} \equiv a_0^{(q-1)/2} \pmod{\frakp} \]
using repeated squaring, where $q=\#k_v$.

To solve Problem \probsf{HasPoint}, by Proposition \ref{padicdet} we have two cases.  In the first case, where one value of the square symbol is equal to $1$, we reduce to Problem \probsf{SquareRoot} over $F_v^{\textup{alg}}$ which we can solve by the above.  Otherwise, if all three symbols in Proposition \ref{padicdet} are $-1$, then also by Hensel's lemma, Problem \probsf{ExhibitPoint} over $F_v^{\textup{alg}}$ is reducible to Problem \probsf{ExhibitPoint} over $k_v$, which was discussed at the end of the previous section.

If we restrict our input to a global field $F$, then a runtime analysis of the above method yields the following.

\begin{prop} \label{Hilbertsymbolodd}
Let $F$ be a number field and let $v$ be an odd place of $F$.  Then there exists a deterministic polynomial-time algorithm to evaluate the Hilbert symbol $(a,b)_v$ for $a,b \in F^\times$.
\end{prop}

\begin{rmk}
By \emph{Hilbert reciprocity}, we have
\begin{equation} \label{Hilbertrecip}
 \prod_v (a,b)_v = 1
\end{equation}
whenever $F$ is a global field and $a,b \in F^\times$.  Consequently, if one can compute all but one local Hilbert symbol $(a,b)_v$, then the final symbol can be recovered from the above relation.  In particular, this means for a number field $F$, if there exists a unique prime above $2$ (e.g.\ when $F=\Q$) then one can evaluate $(a,b)_2$ in this way.  
\end{rmk}

\section{The even local Hilbert symbol}

In this section, we discuss the computation of the local Hilbert symbol for an even place of a number field $F$.  The main result of this section is the following theorem.
 

\begin{thm} \label{evenhilb}
Let $F$ be a number field and let $v$ be a place of $F$.  Then there exists a deterministic polynomial-time algorithm to evaluate the Hilbert symbol $(a,b)_v$ for $a,b \in F^\times$.
\end{thm}

If $v$ is complex, this theorem is trivial; if $v$ is an odd place of $F$ then Theorem \ref{evenhilb} follows from Proposition \ref{Hilbertsymbolodd}.  So suppose that $v$ is an even place of $F$, i.e.\ $\#k_v$ is even.  Let $\Z_F$ be the ring of integers of $F$ and let $\frakp$ be the prime of $\Z_F$ corresponding to $v$.

We first give an algorithm which gives a solution to an integral norm form via a Hensel-like lift.

\begin{alg} \label{evennorm}
Let $\frakp$ an even prime with ramification index $e = \ord_{\frakp} 2$, and let $a,b \in F$ be such that $\ord_\frakp(a)=0$ and $\ord_\frakp(b)=1$.  This algorithm outputs a solution to the congruence
\[ 1-ay^2-bz^2 \equiv 0 \pmod{\frakp^{2e}} \]
with $y,z \in \Z_F/\frakp^{2e}$ and $y \in (\Z_F/\frakp)^{\times}$.

\begin{enumalg}
\item Let $f \in \Z_{\geq 1}$ be the residue class degree of $\frakp$ (so that $\#(\Z_F/\frakp)=2^f$) and let $q=2^f$.  Let $\pi$ be a uniformizer at $\frakp$.

\item Initialize $(y,z) := (1/\sqrt{a},0)$.

\item Let $N := 1-ay^2-bz^2 \in \Z_F/4\Z_F$ and let $t := \ord_\frakp(N)$.  If $t \geq 2e$, return $y,z$.  Otherwise, if $t$ is even, let
\[ y := y+\sqrt{\frac{N}{a\pi^t}} \pi^{t/2} \]
and if $t$ is odd, let
\[ z := z+\sqrt{\frac{N}{b\pi^{t-1}}} \pi^{\lfloor t/2 \rfloor}. \]
Return to Step 3.
\end{enumalg}
\end{alg}

In this algorithm, when we write $\sqrt{u}$ for $u \in (\Z_F/\frakp^{2e})^\times$ we mean any choice of a lift of $\sqrt{u} \in (\Z_F/\frakp)^\times$ to $\Z_F/\frakp^{2e}$.

\begin{proof}[Proof of correctness]
The key calculation in Step $3$ is as follows: if $t$ is even, we make the substitution
\[ 1 - a(y+u\pi^{t/2})^2 - bz^2 = N - 2au\pi^{t/2}y - au^2\pi^t \equiv 0 \pmod{\frakp^{t+1}} \]
and solve for $u$.  Note that since $t < 2e$ we have $\ord_\frakp(2\pi^{t/2})=e+t/2 \geq t+1$; solving we get $u^2 \equiv N/(a\pi^t) \pmod{\frakp}$ as claimed.  The case where $t$ is odd is similar: we have
\begin{align*} 
1-ay^2-b(z+\sqrt{N/b\pi^{t-1}} \pi^{\lfloor t/2 \rfloor})^2 
&= N - 2bz\sqrt{N/b\pi^{t-1}}\pi^{\lfloor t/2 \rfloor} - b(N/b\pi^{t-1})\pi^{t-1}  \\
&\equiv N - N \equiv 0 \pmod{\frakp^{t+1}}
\end{align*}
and the middle term above vanishes modulo $\frakp^{t+1}$ since $t < 2e$ implies $e+1+\lfloor t/2 \rfloor = e+1+(t-1)/2 \geq t+1$.
\end{proof}

\begin{rmk}
Alternatively, we can compute a solution modulo $2$ directly.  The map
\begin{align*} 
(\Z_F/\frakp^e)^2 &\to \Z_F/2\Z_F \\
(y,z) &\mapsto 1-ay^q-bz^q
\end{align*}
is $\Z_F/\frakp \cong \F_q$-linear since $2 \equiv 0 \pmod{\frakp^e}$.  Let $(y_0,z_0)$ be in the kernel of this map.  Letting $(x,y,z) := (1,y_0^{q/2},z_0^{q/2})$, we see $1-ay^2-bz^2 \equiv 0 \pmod{2}$.
\end{rmk}

\begin{rmk}
This is better than the algorithm provided in Simon's thesis \cite{Simonthesis} because we do not need to make a brute force search, which might not run in polynomial time.
\end{rmk}

We reduce to the above Hensel lift by the following algorithm.  

\begin{alg} \label{valuationgame}
Let $\frakp$ an even prime with ramification index $e = \ord_{\frakp} 2$ and let $a,b \in F^\times$ be such that $v(a)=0$ and $v(b) \in \{0,1\}$.  This algorithm outputs $y,z,w \in \Z_F/\frakp^{2e}$ such that 
\[ 1-ay^2-bz^2+abw^2 \equiv 0 \pmod{\frakp^{2e}} \]
and $y \in (\Z_F/\frakp)^{\times}$.  Let $\pi$ be a uniformizer for $\frakp$.

\begin{enumalg}
\item If $v(b)=1$, return the output $(y,z,0)$ of Algorithm \ref{evennorm} with input $a,b$.

\item Suppose $a \in (\Z_F/\frakp^e\Z_F)^{\times 2}$ and $b \in (\Z_F/\frakp^e\Z_F)^{\times 2}$.  Let $(a_0)^2 a \equiv 1 \pmod{\frakp^e}$ and $(b_0)^2 b \equiv 1 \pmod{\frakp^e}$.  Return 
\[ y :=a_0,\ z :=b_0,\ w :=a_0b_0. \]

\item Swap $a,b$ if necessary so that $a \in (\Z_F/\frakp^e\Z_F)^\times \setminus (\Z_F/\frakp^e\Z_F)^{\times 2}$.  Let $t$ be the largest integer such that $a \in (\Z_F/\frakp^t)^{\times 2}$ but $a \not\in (\Z_F/\frakp^e)^{\times 2}$.  Then $t$ is odd; write $a = a_0^2 + \pi^t a_t$ with $a_0,a_t \in \Z_F$.  Let $y,z$ be the output of Algorithm \ref{evennorm} with input $a' := a$, $b' := -\pi a_t/b$.  Return 
\[ y':=\frac{1}{a_0},\ z':=\frac{\pi^{\lfloor t/2 \rfloor}}{a_0 z},\ w':=\frac{y\pi^{\lfloor t/2 \rfloor}}{a_0 z} \] 
(reswapping if necessary).
\end{enumalg}
\end{alg}

\begin{proof}[Proof of correctness]
In Step 2, writing $aa_0^2 = 1+2a'$ and $bb_0^2 = 1+2b'$ with $a',b' \in \Z_F$ we indeed have
\[ 1-a(a_0)^2-b(b_0)^2 + ab(a_0 b_0)^2 = 1 - (1+2a') - (1+2b') + (1+2a')(1+2b') \equiv 0 \pmod{\frakp^{2e}} \]
since $4 \in \frakp^{2e}$.

Now we discuss Step 3.  Write $a=a_0+a_1\pi+\dots+a_{e-1}\pi^{e-1}$ with $a_i \in \Z_F/\frakp$.  Then indeed $a \in (\Z_F/\frakp^{e})^{\times 2}$ if and only if and $a_i=0$ for $i$ odd by the freshperson's dream, so in particular $t<e$ is odd.  Now suppose from Algorithm \ref{evennorm} we have
\[ 1-ay^2 + (\pi a_t/b)z^2 \equiv 0 \pmod{\frakp^{2e}}. \]
Note $\ord_\frakp(z) \leq \lfloor t/2 \rfloor = (t-1)/2$ since otherwise $a \in (\Z_F/\frakp^{t+1})^{\times 2}$, a contradiction.  Multiplying by $-b \pi^{t-1}/z^2=-b(\pi^{\lfloor t/2 \rfloor}/z)^2$ gives
\[ -b(\pi^{\lfloor t/2 \rfloor}/z)^2 + ab(y \pi^{\lfloor t/2 \rfloor}/z)^2 - \pi^t a_t \equiv 0 \pmod{\frakp^{2e}} \]
so
\[ a_0^2 - (a_0^2+\pi^t a_t) - b(\pi^{\lfloor t/2 \rfloor}/z)^2 + ab(y\pi^{\lfloor t/2 \rfloor}/z)^2 \equiv 0 \pmod{\frakp^{2e}} \]
so since $a=a_0^2+\pi^t a_t$, dividing by $a_0^2$ we have the result.
\end{proof}

We say that $\pi^{-1} \in F$ is an \emph{inverse uniformizer} for $\frakp$ if $\ord_\frakp(\pi^{-1})=-1$ and $\ord_\frakq(\pi^{-1}) \geq 0$ for all $\frakq \neq \frakp$.

We are now prepared to evaluate the even Hilbert symbol.

\begin{alg} \label{evenhilbalg}
Let $B=\quat{a,b}{F}$ be a quaternion algebra with $a,b \in F^{\times}$, and let $\frakp$ be an even prime of $F$.  This algorithm returns the value of the Hilbert symbol $(a,b)_\frakp$.

\begin{enumalg}
\item Scale $a,b$ if necessary by an element of $\Q^{\times 2} \cap \Z$ so that $a,b \in \Z_F$.  

\item Let $\pi^{-1}$ be an inverse uniformizer for $\frakp$.  Let $a := (\pi^{-1})^{2\lfloor \ord_\frakp(a)/2 \rfloor} a$ and $b := (\pi^{-1})^{2\lfloor \ord_\frakp(b)/2 \rfloor} b$.  If $\ord_\frakp a=\ord_\frakp b = 1$, let $a := (\pi^{-1})^2 (-ab)$.  Swap if necessary so that $\ord_\frakp a = 0$.

\item Call Algorithm \ref{valuationgame}, and let $i' := (1+yi+zj+wij)/2$.  Let $f(T)=T^2-T+\nrd(i')$ be the minimal polynomial of $i'$.  If $f$ has a root modulo $\frakp$, return $1$.

\item Let $j' := (zb)i - (ya)j$ and let $b' := (j')^2$.  If $\ord_v b'$ is even, return $1$, otherwise return $-1$.
\end{enumalg}
\end{alg}

\begin{proof}[Proof of correctness]
If in Step 2 we have a root modulo $\frakp$, then by Hensel's lemma, $f$ has a root $t \in F_\frakp$, hence $t-i'$ is a zero divisor and we return $1$ correctly.  Otherwise, by Lemma \ref{vignlem}, we have $K_\frakp=F_\frakp[i']$ is the unramified field extension of $F_\frakp$.  We compute that $\trd(j')=\trd(i'j')=0$, so $B_\frakp \cong \quat{K_\frakp,b'}{F_\frakp}$ and $B_\frakp$ is split if and only if $\ord_\frakp b'$ is even.  
\end{proof}

Note that the above algorithms run in deterministic polynomial time.

\begin{exm}
Let $F=\Q(u)$ where $u=\sqrt[8]{500}$.  Then $2\Z_F = (2, \sqrt[8]{500})^4=\frakp^4$, so $\Z_{F,\frakp}$ is a ramified extension of $\Z_2$ of residue degree $2$ and ramification degree $e=4$.  Using Algorithm \ref{evenhilbalg}, we compute $(a,b)_\frakp$ where $b=u^2+40$ and $a=u^2+u+1$.

In Step 2, we compute the inverse uniformizer $\pi^{-1} = u^3/10$ satisfying the polynomial $T^8 - 5/4$.  We compute $\ord_\frakp(a) = 0$ and $\ord_\frakp(b)=2$.  So we let $b := (\pi^{-1})^2 b=\frac{1}{5}(2u^6 + 25)$ with now $\ord_\frakp(b)=0$.

In Step 3, we call Algorithm \ref{valuationgame}.  We use the uniformizer $\pi=u$.  We compute that $b \equiv 1 \pmod{\frakp^e}$ so $b \in (\Z_F/\frakp^e\Z_F)^{\times 2}$ but $a \equiv 1+\pi + \pi^2 \pmod{\frakp^e}$.  So we write $a=a_0+\pi^t a_t$ with $a_0=1$ and $a_t=u+1$.  

We then call Algorithm \ref{evennorm} with input $a' := a$ and $b' := -\pi a_t/b$.  We initialize $(y,z) = (1,0)$.  In Step 3 of this algorithm, we have $N := 1-(1+u+u^2)=-(u+u^2)$ with valuation $t := 1$.  We let $z := \sqrt{N/b} = 1$ and return; now $N := 1 - ay^2 - bz^2$ has valuation $t := 9 > 2e$, so we exit the loop with output $y = z = 1$.  

We then exit Algorithm \ref{valuationgame} with $y' := 1/a_0 = 1$, $z' := \pi^{\lfloor t/2 \lfloor}/(a_0 z) = 1$, and $w' := y \pi^{\lfloor t/2 \lfloor}/(a_0 z) = 1$.  We verify that $1 - a(y')^2 - b(z')^2 + ab(w')^2 = 1-a-b+ab \equiv 0 \pmod{4}$.  

Returning to Algorithm \ref{evenhilbalg}, we let $i' := (1+i+j+ij)/2$ and compute $\nrd(i') = 1/10(w^7+10w^2+10w+500) \equiv 0 \pmod{\frakp}$, so $f(T)=T^2-T+\nrd(i')$ has a root modulo $\frakp$, and we return $(a,b)_\frakp=1$.  
\end{exm}

\subsection*{Computing the Jacobi symbol}

An interesting consequence of the above algorithm is that one can evaluate the Jacobi symbol in deterministic polynomial time in certain cases analogous to the way (``reduce and flip'') that one computes this symbol using quadratic reciprocity in the case $F=\Q$.  (See Lenstra \cite{Lenstrasymbol} for an alternative approach which works in greater generality.)

We extend the definition of the Legendre symbol (\ref{Legendredef}) to a symbol $\displaystyle{\legen{a}{\frakb}}$ with $\frakb$ odd by multiplicativity, and we define $\displaystyle{\legen{a}{b}=\legen{a}{b\Z_F}}$.

We write $v \mid 2\infty$ for the set of finite even places and real archimedean places of $F$.

\begin{prop} \label{quadrecip}
Let $a,b \in \Z_F$ satisfy $a\Z_F+b\Z_F=\Z_F$, with $b$ odd, and suppose $a=a_0a_1$ with $a_1$ odd.  Then
\[ \legen{a}{b}\legen{b}{a_1}= \prod_{v \mid 2\infty} (a,b)_v. \]
\end{prop}

\begin{proof}
By Hilbert reciprocity (\ref{Hilbertrecip}), we have
\[ \prod_v (a,b)_v = 1 = \prod_{v \mid 2\infty} (a,b)_v \prod_{\frakp \nmid 2} (a,b)_\frakp. \]
By Lemma \ref{padicdet}, if $\frakp$ is odd and $\ord_\frakp(a)=\ord_\frakp(b)=0$ then $(a,b)_\frakp=1$.  Therefore
\[ \prod_{\frakp \mid a_1 b} (a,b)_\frakp = \prod_{v \mid 2\infty} (a,b)_v. \]
For $\frakp$ odd, if $\ord_\frakp a_1 > 0$ then $\ord_\frakp b=0$ by assumption and thus
\[ (a,b)_\frakp=\legen{b}{\frakp}^{\ord_\frakp a} = \legen{b}{\frakp}^{\ord_\frakp a_1}. \]
Similarly if $\ord_\frakp b > 0$ then $(a,b)_\frakp=\legen{a}{\frakp}^{\ord_\frakp b}$, hence
\[ \prod_{\frakp \mid a_1 b} (a,b)_\frakp = \legen{a}{b}\legen{b}{a_1}. \]
The result follows.
\end{proof}

A \emph{Euclidean function} on $F$ is a map $N:\Z_F \setminus \{0\} \to \Z_{\geq 0}$ such that for all $a,b \in \Z_F$ we have $N(ab)=N(a)N(b)$ and there exists $q,r \in \Z_F$ such that $a=qb+r$ with either $r=0$ or $N(r)<N(b)$.  A Euclidean function is \emph{computable} if given $a,b$, the elements $q,r$ as above are computable.

\begin{alg} \label{jacobi}
Let $F$ be a number field with a computable Euclidean function $N$ and let $a,b \in \Z_F \setminus \{0\}$.  This algorithm returns the Jacobi symbol $\displaystyle{\legen{a}{b}}$.

\begin{enumalg}
\item Initialize $z=1$.
\item If $b\Z_F=\Z_F$, return $z$.  Otherwise, compute $q,r \in \Z_F$ such that $a=qb+r$.  If $r=0$, return $0$.  Let $a := r$.  Write $a=a_0 a_1$ with $a_1 \in \Z_F$ odd.  
\item Multiply $z$ by $\prod_{v \mid 2,\infty} (a,b)_v$, computed using Algorithm \ref{evenhilbalg}.  Return to Step 2, with $(a,b)=(b,a_1)$.
\end{enumalg}
\end{alg}

\begin{proof}[Proof of correctness]
The division algorithm associated to $N$ implies that $\Z_F$ has unique factorization, so we can indeed write $a=a_0 a_1$ with $a_1$ odd.  The algorithm terminates because in Step $4$ we have $N(a_1) \leq N(a)=N(r)<N(b)$.
\end{proof}

\begin{rmk}
For any fixed $F$, one can precompute a table of the values $(a,b)_{\frakp}$ for $a,b$ in appropriate residue classes modulo an even prime $\frakp$; this is what is usually done for $F=\Q$, for example.
\end{rmk} 

\subsection*{Relationship to conics}

In view of the results in Section 4, we now relate the above algorithms to the geometric problem of rational points on conics.

\begin{thm}[Hasse-Minkowski]
A quaternion algebra $B$ has $B \cong \M_2(F)$ if and only if $B$ is unramified at all places.  
\end{thm}

Equivalently, a conic $C$ has $C(F) \neq \emptyset$ if and only if $C(F_v) \neq \emptyset$ for all places $v$ of $F$.  For a proof of the Hasse-Minkowski Theorem, see Lam \cite{Lam}, O'Meara \cite{OMeara}, or Vign\'eras {\cite[\S III.3.1]{Vigneras}}

\begin{prop}
Problem \textup{(\textsf{IsMatrixRing})} is deterministic polynomial-time reducible to the problem of factoring ideals in $\Z_F$.
\end{prop}

\begin{proof}
Given a quaternion algebra $B=\quat{a,b}{F}$, we have $B_v \cong \M_2(F_v)$ for all $v \nmid 2ab\infty$, and by factoring by the above algorithms for each $v \mid 2ab\infty$ we check if $B_v \cong \M_2(F_v)$ by computing the Hilbert symbol $(a,b)_v$ in deterministic polynomial time.
\end{proof}

\section{Maximal orders}

In this section, we consider some integral versions (for orders) of the above algorithms relating quadratic forms and quaternion algebras.  Our main result relates identifying the matrix ring to computing a maximal order.  Throughout this section, let $F$ be a number field, let $\Z_F$ be its ring of integers, and let $\calO$ be a ($\Z_F$-)order in a quaternion algebra $B$ over $F$.  For further reading, see Reiner \cite{Reiner} or Vign\'eras \cite{Vigneras}.

\subsection*{Computing maximal orders, generally}

There exists a deterministic algorithm to compute the ring of integers $\Z_F$ (see Cohen \cite[\S 6.1]{Cohen}, \cite[Algorithm 2.4.9]{Cohen2}): in fact, computing $\Z_F$ is deterministic polynomial-time equivalent to the problem of finding the largest square divisor of a positive integer \cite{Chistov,Lenstra}; no polynomial-time algorithm is known for this problem (though see work of Buchmann and Lenstra \cite{BuchLenstra} for a way of ``approximating'' $\Z_F$).

\begin{exm}
If $F=\Q(\sqrt{D})$, then $R=\Z \oplus \Z(d+\sqrt{d})/2$ where $D=df^2$ and $f^2$ is the largest square divisor of $D$ subject to the requirement that $d \equiv 0,1 \pmod{4}$.
\end{exm}

We consider in this section the noncommutative analogues of this problem.  We have the following general result due to Ivanyos and R\'onyai \cite[Theorem 5.3]{Ronyai3}, which was rediscovered by Nebe and Steel \cite{NebeSteel}; see also Friedrichs \cite{Friedrichs}.

\begin{thm} \label{maxorderthm}
There exists an explicit algorithm which, given a semisimple $F$-algebra $B$, computes a maximal order $\calO \subset B$.  This algorithm runs in deterministic polynomial time given oracles for the problems of factoring integers and factoring polynomials over finite fields.  
\end{thm}

At present, it is not known if there exist deterministic polynomial-time algorithms to solve either of these latter two problems.  Indeed, we have already noted that computing a maximal order in $F$ is as hard as computing the largest squarefree divisor of a positive integer; therefore, it is reasonable to expect that the problem for a noncommutative algebra $B$ is no less complicated.  (See a more precise characterization of this complexity at the end of this section.)

We do not discuss the algorithm exhibited in Theorem \ref{maxorderthm}; rather, we consider the special case of quaternion algebras, and by manipulations with quadratic forms we obtain a simpler algorithm.

\subsection*{Discriminants}

We begin by analyzing the following problem.  

\begin{prob}[\textsf{IsMaximal}] \label{ismaximal}
Given an order $\calO \subset B$, determine if $\calO$ is a maximal order.
\end{prob}

This problem has a very simple solution as follows.  The \emph{discriminant} $\frakD(B)$ of $B$ is the ideal equal to the product of all primes of $\Z_F$ where $B$ is ramified:
\[ \frakD(B) = \prod_{\frakp \text{ ramified}} \frakp. \]
On the other hand, the \emph{discriminant} $\disc(\calO)$ of an order $\calO \subset B$ is the ideal generated by the set
\[ \{ \det( \trd(x_ix_j) )_{i,j=1,\dots,4} : x_1,\dots,x_4 \in \calO \}. \]
The discriminant $\disc(\calO)$ is the square of an ideal in $\Z_F$, and the \emph{reduced discriminant} $\frakd(\calO)$ of $\calO$ is the ideal satisfying $\frakd(\calO)^2=\disc(\calO)$.

Given a pseudobasis $(\fraka_i, x_i)$ for $\calO$ we have
\[ \disc(\calO)=(\fraka_1 \cdots \fraka_4)^2 \det( \trd(x_ix_j) )_{i,j=1,\dots,4}. \]

\begin{rmk}
Although we will not use this in the sequel, the reduced discriminant can in fact be computed more simply: if $\calO=\Z_F \oplus \fraka i \oplus \frakb j \oplus \frakc k$ then
\[ \frakd(\calO)=\fraka\frakb\frakc \trd((ij-ji)\overline{k}). \]
\end{rmk}

\begin{lem} \label{disccharmaxorder}
An order $\calO \subset B$ is maximal if and only if $\frakd(\calO)=\frakD(B)$.
\end{lem}

\begin{proof}
We give only a sketch of the proof.  For a prime $\frakp$ of $\Z_F$, let $\Z_{F,\frakp}$ be the completion of $\Z_F$ at $\frakp$ and $F_\frakp$ the completion of $F$ at $\frakp$; write $\calO_\frakp = \calO \otimes_{\Z_F} \Z_{F,\frakp}$ and similarly $B_\frakp = B \otimes_{F} F_\frakp$.  

We have $\frakd(\calO)=\frakD(B)$ if and only if $\frakd(\calO)_\frakp=\frakd(\calO_\frakp)=\frakD(B_\frakp)=\frakD(B)_\frakp$ for all primes $\frakp$, and the order $\calO$ is maximal if and only if $\calO_\frakp$ is maximal for every prime $\frakp$ of $\Z_F$ (see \cite[11.2]{Reiner}).  So it suffices to note that if $\frakp$ is unramified then any maximal order of $B_\frakp$ has discriminant $\Z_{F,\frakp}$ and if $\frakp$ is ramified then the unique maximal order of $B_\frakp$ has reduced discriminant $\frakp \Z_{F,\frakp}$ \cite[Theorem 14.9]{Reiner}.
\end{proof}

Putting these together with the computation of the local Hilbert symbol, we have shown that one can solve Problem \probsf{IsMaximal} in deterministic polynomial time given an oracle to factor integers and polynomials over finite fields, since this allows the factorization of the discriminant $\frakD(B)$ \cite[Proposition 6.2.8, Algorithm 6.2.9]{Cohen}; note that this need only be done once for a quaternion algebra $B$.

\subsection*{Computing maximal orders}

We now turn to the problem of computing a maximal order in a quaternion algebra.

\begin{prob}[\textsf{AlgebraMaxOrder}] \label{findmaxorder}
Given a quaternion algebra $B$ over $F$, compute a maximal order $\calO \subset B$.
\end{prob}

A more general problem is as follows.

\begin{prob}[\textsf{MaxOrder}] \label{findmaxordercontaining}
Given an order $\Lambda \subset B$ in a quaternion algebra $B$ over $F$, compute a maximal order $\calO \supset \Lambda$.
\end{prob}

One immediately reduces from the former to the latter by exhibiting any order in $B$, as follows.  (First, we compute $\Z_F$ as above; this can be considered a precomputation step if $F$ is fixed.)  If $B=\quat{a,b}{F}$, we may scale $a,b$ by a nonzero square integer so that $a,b \in \Z_F$, and then 
\begin{equation} \label{Lambdayeah}
\Lambda = \Z_F \oplus \Z_F i \oplus \Z_F j \oplus \Z_F ij
\end{equation} 
is an order, where $i,j$ are the standard generators for $B$.

An order $\calO$ is \emph{$\frakp$-maximal} for a prime $\frakp$ if $\calO_\frakp=\calO \otimes_{\Z_F} \Z_{F,\frakp}$ is maximal (as an $\Z_{F,\frakp}$-order).  Note that if $\ord_\frakp(\frakd(\calO_\frakp))=0$ then necessarily $\calO$ is $\frakp$-maximal.  To solve Problem \probsf{MaxOrder}, we recursively compute a $\frakp$-maximal order for every prime $\frakp \mid \frakd(\calO)$, proceeding in two steps.  

We say an order $\calO$ is \emph{$\frakp$-saturated} if $\nrd|_{\calO_\frakp}$ has a normalized basis $1,i,j,k$ (see Proposition \ref{qformdiag}) such that each atomic block has valuation at most $1$; we then say that $1,i,j,k$ is a \emph{$\frakp$-saturated} basis for $\calO$.

We compute a $\frakp$-saturated order in the following straightforward way.  Recall that $\pi^{-1} \in F$ is an \emph{inverse uniformizer} for $\frakp$ if $\ord_\frakp(\pi^{-1})=-1$ and $\ord_\frakq(\pi^{-1}) \geq 0$ for all $\frakq \neq \frakp$.

\begin{alg} \label{computesaturatedorder}
Let 
\[ \Lambda = \Z_F \oplus \fraka i \oplus \frakb j \oplus \frakc k \subset B \] 
be an order and let $\frakp$ be prime.  This algorithm computes a $\frakp$-saturated order $\calO \supset \Lambda$ and a $\frakp$-saturated basis for $\calO$.

\begin{enumalg}
\item Choose $d \in \fraka$ such that $\ord_\frakp(d) = \ord_\frakp(\fraka)$ and let $i := di$; compute similarly with $j$, $k$.  Let $\calO := \Lambda$.

\item Run Algorithm \ref{normalizequadform} over the localization of $\Z_F$ at $\frakp$ with input the quadratic form $\nrd|_{\calO}$ and the basis $1,i,j,k$; let $1,i^*,j^*,k^*$ be the output.  Let $c \in \Z_F$ be such that $\ord_\frakp c = 0$ and such that $ci^* \in \calO$, and let $i := ci^*$; compute similarly with $j$, $k$.   

\item Let $\pi^{-1}$ be an inverse uniformizer for $\frakp$.  For each atomic form $Q$ in $\nrd_{\calO}$, let $e$ be the valuation of $Q$, and multiply each basis element in $Q$ by $(\pi^{-1})^{\lfloor e/2 \rfloor}$.  Return $\calO := \Lambda + (\Z_F i \oplus \Z_F j \oplus \Z_F k)$ and the basis $1,i,j,k$.
\end{enumalg}
\end{alg}

\begin{proof}[Proof of correctness]
In Step 3, we are asserting that the output of Algorithm \ref{normalizequadform} leaves $1$ as the first basis element.  Indeed, we note that $\ord_\frakp \trd(j) \leq \ord_\frakp \trd(i(ij))$ since $\trd(i(ij))=\trd(i)^2-\trd(j)\nrd(i)$ and similarly $\ord_\frakp \trd(i) \leq \ord_\frakp \trd((ij)j)$.  

Let $1,i,j,k$ be the basis computed in Step 3.  By definition, this basis is $\frakp$-saturated; we need to show that $\calO$ is indeed an order.  But $\calO$ is an order if and only if $\calO_\frakq$ is an order for all primes $\frakq$, and we have $\calO_\frakq=\Lambda_\frakq$ for all primes $\frakq \neq \frakp$.  

For any $x,y \in B$ we have $xy+yx=\trd(y)x+\trd(x)y-T(x,y)$, so if $\calO$ is an order then $\calO + \Z_F x$ is multiplicatively closed if and only if $T(x,y) \in \Z_F$ for all $y \in \calO$.  We have $T(x,y)=0$ if $x,y$ are orthogonal, and if $x,y$ are a basis for an atomic block $Q$ then by definition the valuation of $T(x,y)$ is at least the valuation of $Q$ and so we can multiply each by $(\pi^{-1})^{\lfloor e/2 \rfloor}$, preserving integrality.
\end{proof}

After $\frakp$-saturating, one can compute a maximal order as follows.

\begin{alg} \label{computepmaxorder}
Let $\Lambda$ be an order and let $\frakp$ be prime.  This algorithm computes a $\frakp$-maximal order $\calO \supset \Lambda$.

\begin{enumalg}
\item Compute a $\frakp$-saturated order $\calO \supset \Lambda$ and let $1,i,j,k$ be a $\frakp$-saturated basis for $\calO$.  Let $\pi^{-1}$ be an inverse uniformizer for $\frakp$.  

\item Suppose $\frakp$ is odd.  Swap $i$ for $j$ or $k$ if necessary so that $a := i^2$ has $\ord_\frakp(a)=0$.  Let $b := j^2$.  If $\ord_\frakp b = 0$, return $\calO$.
Otherwise, if $\ord_\frakp b = 1$ and $(a/\frakp)=1$, solve
\[ x^2 \equiv a \pmod{\frakp} \]
for $x \in \Z_F/\frakp$.  Adjoin the element $\pi^{-1} (x-i)j$ to $\calO$, and return $\calO$.

\item Otherwise, $\frakp$ is even.  Let $t := \trd(i)$, let $a := -\nrd(i)$, and let $b := j^2$.
\begin{enumalgalph}
\item Suppose $\ord_\frakp t=0$.  If $\ord_\frakp b = 0$, return $\calO$.  If $\ord_\frakp b = 1$ and $T^2- t T + a=0$ has a root $x$ modulo $\frakp$, and return $\calO + \Z_F \pi^{-1} (x-i)j$.  

\item Suppose $\ord_\frakp \trd(i)>0$.  Let $y,z,w$ be the output of Algoritm \ref{valuationgame} with input $a,b$.  Let 
\[ i' := (\pi^{-1})^e(1+yi+zj+wij). \]
Adjoin $i'$ to $\calO$, and return to Step 1.
\end{enumalgalph}
\end{enumalg}
\end{alg}

\begin{proof}[Proof of correctness]
At every step in the algorithm, for each prime $\frakq \neq \frakp$ the order $\calO_\frakq$ does not change, so we need only verify that $\calO_\frakp$ is indeed a maximal order.

In Step 2, we have that $b$ is a uniformizer for $\frakp$, that $\frakd(\calO_\frakp)=4ab\Z_{F,\frakp}$.  If $\ord_\frakp(b)=0$ then $\ord_\frakp \frakd(\calO_\frakp) = 0$ so $\calO$ is indeed maximal.  Otherwise, we have $\frakd(\calO_\frakp)=\frakp$ and $B_\frakp \cong \quat{K_\frakp,b}{F_\frakp}$ where $K_\frakp=F_\frakp[i]$.  We conclude that $B_\frakp$ is a division ring (and hence $\calO_\frakp$ is maximal) if and only if $(a/\frakp)=-1$.  If $(a/\frakp)=1$ and $j'=\pi^{-1}(x-i)j$, then $1,i,j',ij'$ form the $\Z_{F,\frakp}$-basis for a maximal order, since $(j')^2=(\pi^{-1})^2(x^2-a)b \in \Z_{F,\frakp}$ and $j'i = -ij'$.  

In Step 3, first note that $ij$ is also orthogonal to $1,i$: we have $i$ orthogonal to $j$ so $\trd(ij)=0$ so $ij$ is orthogonal to $1$, and similarly $\trd(ij\overline{i})=\trd(\nrd(i)j)=0$.  In particular, we have $B_\frakp=\quat{K_\frakp,b}{F_\frakp}$ where $K_\frakp=F_\frakp[i]$.  By a comparison of discriminants, using the fact that the basis is normalized, we see that $1,i,j,ij$ is a $\frakp$-saturated basis for $\calO$ as well, so without loss of generality we may take $k=ij$.  

Suppose first that $\ord_\frakp \trd(i)=0$, so we are in Step 3a.  If $\ord_\frakp b=0$, then $\ord_\frakp \frakd(\calO_\frakp)=0$ so $\calO_\frakp$ is maximal.  If $\ord_\frakp b>0$, then since the basis is $\frakp$-saturated we have $\ord_\frakp b = 1$.  Thus as in the case for $\frakp$ odd, we have $B_\frakp$ is a division ring if and only if $K_\frakp$ is not a field, and as above the adjoining the element $\pi^{-1}(x-i)j$ yields a maximal order.

So suppose we are in Step 3b, so $\ord_\frakp \trd(i)>0$.  Since $1,i,j,k$ is normalized, we have $\ord_\frakp \trd(i)=\ord_\frakp T(1,i) \leq \ord_\frakp T(j,k)$.
Adjoining $i'$ to $\calO$ gives a $\Z_{F,\frakp}$-module with basis $1,i',j,i'j$ since $y \in (\Z_F/\frakp)^\times$; adjoining $j'$ gives a module with basis $1,i',j',i'j'$ for the same reason.  We verify that $\calO_\frakp$ after these steps is indeed an order: we have $\trd(i')=2(\pi^{-1})^e \in \Z_{F,\frakp}$ and $\nrd(i')=(\pi^{-1})^{2e}(1-ay^2-bz^2+abw^2) \in \Z_{F,\frakp}$ by construction, so at least $\Z_{F,\frakp}[i] = \Z_{F,\frakp} \oplus \Z_{F,\frakp} i$ is a ring.  Similarly we have $(j')^2 = b' \in \Z_{F,\frakp}$.  Finally, we have $\trd(i'i)=2(\pi^{-1})^e ya$ and $\trd(i'j)=2(\pi^{-1})^e zb$, so it follows that $\trd(i'j')=0$, and hence $j'i'=-\overline{i'}j'=-i'j'-\trd(i')j'$, so indeed we have an order.  
\end{proof}

\begin{rmk}
One must really treat the even and odd prime cases separately.  Consider, for example, $F=\Q$, and the quaternion algebra $B=\quat{-3,5}{\Q}$.  Then we have the maximal orders $\Z[(1+i)/2] \subset \Q(i) \cong \Q(\sqrt{-3})$ and $\Z[(1+j)/2] \subset \Q(j) \cong \Q(\sqrt{5})$, but we find that
\[ \left(\frac{1+j}{2}\right)\left(\frac{1+i}{2}\right)=\left(\frac{1-i}{2}\right)\left(\frac{1+j}{2}\right)+
\frac{ij}{2}, \]
which is not integral (since $ij/2$ has norm $15/4$).
\end{rmk}

\begin{rmk}
In the proof of correctness for Algorithm \ref{computepmaxorder}, in each case where $\frakp$ is ramified in $B$ we have in fact written $B_\frakp \cong \quat{K_\frakp,\pi}{F_\frakp}$ where $K_\frakp$ is the unramified extension of $F_\frakp$.  The reader will note the similarity between this algorithm and the algorithm to compute the Hilbert symbol: the former extends the latter by taking a witness for the fact that the algebra is split, namely a zerodivisor modulo $\frakp$, and uses this to compute a larger order (giving rise therefore to the matrix ring).  
\end{rmk}

Combining these two algorithms, we have the following immediate corollary.

\begin{cor}
There exists an algorithm to solve \textup{(\textsf{ExhibitMatrixRing})} for orders over $\Z_{F,\frakp}$.
\end{cor}

(We recall the discussion in Section 4 for the representation of local fields and rings.)  In other words, if $\calO \subset B$ is an order in a quaternion algebra $B$ over a number field $F$ and $\frakp$ is prime of $\Z_F$ which is unramified in $B$, then there exists an algorithm to compute an explicit embedding $\calO \hookrightarrow \M_2(\calO_\frakp)$.

Putting these two algorithms together, we have proved the following theorem.

\begin{thm} \label{detredfactoridealmaxorder}
Problem \textup{(\textsf{MaxOrder})} is deterministic polynomial-time reducible to the problem of factoring ideals in $\Z_F$.
\end{thm}

\begin{proof}
Given any order $\Lambda$, we factor its discriminant $\frakd(\Lambda)$, and for each prime $\frakp \mid \frakd(\Lambda)$, we compute a $\frakp$-saturated order containing $\Lambda$ from Algorithm \ref{computesaturatedorder} and a $\frakp$-maximal order $\calO$ containing it using Algorithm \ref{computepmaxorder}.
\end{proof}

\subsection*{Complexity analysis}

Given Theorem \ref{detredfactoridealmaxorder}, we prove the following result which characterizes the abstract complexity class of this problem, following a hint of Ronyai \cite[\S 6]{Ronyai00}.

\begin{thm} \label{equivmaxorder}
Problem \textup{(\textsf{AlgebraMaxOrder})} for any fixed number field $F$ is probabilistic polynomial-time equivalent to the problem of factoring integers.
\end{thm}


To prove the theorem, we will use two lemmas.  The first lemma is a standard fact.

\begin{lem} \label{factorideals}
The problem of factoring integral ideals $\fraka$ of an arbitrary number field is probabilistic polynomial-time equivalent to the problem of factoring integers.
\end{lem}

\begin{proof}
Suppose $\fraka$ is an integral ideal of $F$.  After factoring the absolute discriminant $d_F$ of $F$, we can in deterministic polynomial time compute the ring of integers $\Z_F$ of $F$ as above.  Now let $\fraka$ be an ideal with norm $\N(\fraka)=a$.  After we factor $a$, for each prime $p \mid a$, we decompose $p\Z_F=\prod_i \frakp_i^{e_i}$ into primes by a probabilistic polynomial time algorithm due to Buchmann and Lenstra \cite[Algorithm 6.2.9]{Cohen}: this algorithm uses a probabilistic algorithm to factor polynomials over a finite field, such as the Cantor-Zassenhaus algorithm; see von zur Gathen and Gerhard \cite[Theorem 14.14]{GathenGerhard} or Cohen \cite[\S 3.4]{Cohen}.  (In fact, for our applications, it suffices to have an algorithm to compute a square root in a finite field, for which we may use the algorithm of Tonelli and Shanks \cite[\S 1.5.1]{Cohen}.)

From this list of primes we easily obtain the factorization of $\fraka$.  Conversely, if one has an algorithm to factor ideals, then one may factor $a\Z_F$ into primes and computing norms we recover the prime factorization of $a$ over $\Z$.
\end{proof}

\begin{rmk}
Deterministically, already the problem of finding a nonsquare modulo a prime $p$ is difficult; one unconditional result known is that the smallest quadratic nonresidue of a prime $p$ is of size exponential in $\log p$; under condition of a generalized Riemann hypothesis, one can find a quadratic nonresidue which is of polynomial size in $\log p$.
\end{rmk}

We will also make use of one other lemma.  

\begin{lem} \label{crteff}
Let $\fraka$ be an ideal of $\Z_F$ which is odd, not a square, and not a prime power.  Let
\[ S=\left\{b \in (\Z_F/\fraka)^\times: \text{there exist $\frakp^e,\frakq^f \parallel \fraka$ with $\legen{b}{\frakp}^e = -1$ and $\legen{b}{\frakq}^f = 1$} \right\}. \]
Then $\#S \geq \displaystyle{\frac{1}{2}}\#(\Z_F/\fraka)^\times$.
\end{lem}

\begin{proof}
For an ideal $\frakb$, let $\Phi(\frakb)=\#(\Z_F/\frakb)^\times$.  First consider the case where $\fraka=\frakp^e \frakq^f$ is the product of two prime powers.  Without loss of generality, we may assume $e$ is odd.  If $f$ is even, then $b \in S$ if and only if $(b/\frakp)=-1$, so $\#S = \Phi(\frakp^e)/2 \cdot \Phi(\frakq^f)= \Phi(\fraka)/2$.  If $f$ is odd, then $\# S = 2 (\Phi(\frakp^e)/2)(\Phi(\frakq^f)/2) = \Phi(\fraka)/2$.  

To conclude, write $\fraka=\frakp^e \frakq^f \frakb$ with $\frakb$ coprime to $\frakp\frakq$ and $e$ odd.  Then by the preceding paragraph $\# S \geq (1/2) \Phi(\frakp^e \frakq^f) \Phi(\frakb) = \Phi(\fraka)/2$.  
\end{proof}

\begin{proof}[Proof of Theorem \textup{\ref{equivmaxorder}}] 
Since one can factor ideals in probabilistic polynomial time given an algorithm to factor integers by Lemma \ref{factorideals}, we may compute a maximal order as in the previous section as the resulting computations run in (deterministic) polynomial time.

Now we prove the converse.  Suppose we have an algorithm to solve Problem (\textsf{AlgebraMaxOrder}).  Let $a \in \Z_{>0}$ be the integer to be factored, which we may assume without loss of generality is odd, not a prime power, and not a square.  We can in constant time (for fixed $F$) factor the absolute discriminant $d_F$, so we may also assume $\gcd(a,d_F)=1$.  It follows that the ideal $a\Z_F$ is also odd, not a prime power, and not a square.  

We compute a random $b \in \Z_F/a\Z_F$ with $b \neq 0$.  Since $\N(a\Z_F)=a^d$ where $d=[F:\Q]$, if $\N(b\Z_F)$ is not a power of $a$ then dividing $\gcd(a^d,\N(b))$ by powers of $a$ we obtain a factor of $a$.  Otherwise, $\fraka=a\Z_F+b\Z_F$ is a proper divisor of $a\Z_F$, and we repeat, computing a random $b \in \Z_F/\fraka$---in at most $d$ steps, we will either factor $a$ or find an element $b$ such that $a\Z_F + b\Z_F=\Z_F$.  Note $d$ depends only on $F$ and not on $B$, so we find such a $b$ in probabilistic polynomial time.  

By Lemma \ref{crteff}, we can find in probabilistic polynomial time $b \in (\Z_F/a\Z_F)^\times$ such that $\frakp^e,\frakq^f \parallel a$ with $(b/\frakp)^e=-1$ and $(b/\frakq)^f=1$, say.  Let $B=\quat{a,b}{F}$.  By hypothesis, calling an algorithm to solve (\textsf{AlgebraMaxOrder}) we may compute a maximal order $\calO \subset B$.  

We claim that $\frakp \mid \frakd(\calO)$ but $\frakq \nmid \frakd(\calO)$.  Assuming this, we have that $\gcd(\N(\frakd(\calO)),a)$ is a proper factor of $a$, and the proof is complete.

First we prove that $\frakp \mid \frakd(\calO)$.  Since $\frakp$ is prime to $d_F$, we know that $\frakp$ is unramified in $F$, and since $\frakp^e \parallel a\Z_F$ with $e$ odd, the extension $F(\sqrt{a})/F$ is ramified at $\frakp$.  Since $(b/\frakp)=-1$, by Corollary \ref{padicdet}, the algebra $B$ is ramified at $\frakp$.  Therefore by Lemma \ref{disccharmaxorder}, $\frakp$ divides the discriminant $\frakd(\calO)$. 

Now we show that $\frakq \nmid \frakd(\calO)$.  If $f$ is even, since $\frakq^f \parallel a\Z_F$, we have that $F(\sqrt{a})/F$ is unramified at $\frakq$; since also $(b/\frakq) \neq 0$, by the same corollary, $B$ is unramified at $\frakq$.  And if $f$ is odd, then since $(b/\frakq)^f=1$ we must have $(b/\frakq)=1$, and again by the corollary it follows that $B$ is unramified.  
\end{proof}

\subsection*{Relationship to conics}

We return once again to the theme of rational points on conics.

We have seen that given an algorithm to factor integers, one can solve both problems (\textsf{IsMatrixRing}), or equivalently (\textsf{HasPoint}), over a number field $F$ in probabilistic polynomial time by factoring the discriminant and computing Hilbert symbols.  We have also seen that (\textsf{AlgebraMaxOrder}) over a number field $F$ is probabilistic polynomial time equivalent to the problem of factoring integers.  

We are left to consider (\textsf{ExhibitMatrixRing}), or equivalently (\textsf{ExhibitPoint}).  In the special case where $F=\Q$, one shows that again they are reducible to the problem of integer factorization.  

\begin{thm}[Cremona-Rusin \cite{CremRusin}, Ivanyos-Sz\'ant\'o \cite{RonyaiLLL}, Simon \cite{Simonsolve}]
There exists an explicit algorithm to solve \textup{(\textsf{ExhibitPoint})} over $\Q$ which runs in deterministic polynomial time given an oracle to factor integers.
\end{thm}

From our point of view, the algorithm(s) described in the above theorem can be rephrased in the following way: there exists an explicit algorithm which, given a order $\calO$ over $\Z$ of discriminant $1$ which is split at $\infty$, computes a zerodivisor $x \in \calO$.  This algorithm proceeds by computing a reduced basis of $\calO$ with respect to the reduced norm $\nrd$, a kind of indefinite LLL-algorithm.  

\begin{ques}
Does there exist an algorithm which, given an order $\calO$ over $\Z_F$ of discriminant $1$ which is split at all real places of $F$, computes a zerodivisor $x \in \calO$?
\end{ques}

One possible approach to this conjecture, then, is to provide an indefinite LLL algorithm over $F$ in the special case of $\Z_F$-module of rank $4$ and discriminant $1$.  Perhaps one can prove this at least in the case where $\Z_F$ is computably Euclidean?  

We discuss the computational complexity of problem (\textsf{IsMatrixRing}) over $\Q$ in the next section (and relate this to the problem of factoring integers).  From the discussion above, it seems reasonable to conjecture the following.

\begin{conj}
Problem \textup{(\textsf{ExhibitPoint})} over $\Q$ is (probabilistic) polynomial-time equivalent to the problem of factoring integers.  
\end{conj}

Having treated the case of number fields in some detail, we note that over more general fields, the literature is much less complete.  

\begin{ques}
For which computable fields $F$ is there an effective algorithm to solve Problems \probsf{HasPoint} and \probsf{ExhibitPoint}?
\end{ques}

For example, one may ask for which fields $F$ is there an effective version of the Hasse-Minkowski theorem?  Of course, if one can solve \probsf{HasPoint}, then given a conic which is known to have a solution one can always simply enumerate the points of $\PP^2(F)$ until a solution is found.

\section{Residuosity}

In this final section, we return to Problem (\textsf{IsMatrixRing}) and characterize its computational complexity.  Let $F$ be a number field with ring of integers $\Z_F$.  

For a nonzero ideal $\frakb$ of $\Z_F$, let $\sqrad(\frakb)$ be the product of the prime ideals $\frakp$ dividing $\frakb$ to odd exponent, or equivalently the quotient of $\frakb$ by the largest square ideal dividing $\frakb$.

\begin{prob*}[\textsf{QuadraticResiduosity}]
Given an odd ideal $\frakb$ and $a \in \Z_F$, determine if $a \in (\Z_F/\sqrad(\frakb))^{\times 2}$, i.e., determine if $a$ is a quadratic residue modulo $\sqrad(\frakb)$.
\end{prob*}

Problem (\textsf{QuadraticResiduosity}) reduces to the more familiar problem of quadratic residuosity when $\frakb$ is a squarefree ideal, namely, to determine if $a \in (\Z_F/\frakb)^{\times 2}$.  If $\frakb=\frakp$ is a prime ideal, one has $a \in (\Z_F/\frakp)^{\times 2}$ if and only if $(a/\frakp)=1$, and this Legendre symbol can be evaluated in deterministic polynomial time (as discussed above, by repeated squaring).  In general, for $\frakb$ squarefree, we have $a \in (\Z_F/\frakb)^{\times 2}$ if and only if $a \in (\Z_F/\frakp)^{\times 2}$ for all primes $\frakp \mid \frakb$.  In particular, by this reduction if one can factor $\frakb$, one can solve Problem (\textsf{QuadraticResiduosity}).  It is a terrific open problem in number theory to determine if the converse holds, even for the case $F=\Q$ and $\frakb$ generated by $pq$ with $p,q$ distinct primes.

We first relate the problems (\textsf{IsMatrixRing}) and (\textsf{QuadraticResiduosity}) as follows.

\begin{prop}
Problem \textup{(\textsf{IsMatrixRing})} over $F$ is deterministic polynomial-time reducible to Problem \textup{(\textsf{QuadraticResiduosity})} over $F$.
\end{prop}

\begin{proof}
Let $B=\quat{a,b}{F}$ be a quaternion algebra over $F$.  Scaling $a,b$ by an integer square, we may assume $a,b \in \Z_F$.  Recall that $B \cong \M_2(F)$ if and only if every place $v$ of $F$ is unramified in $B$, i.e., if $(a,b)_v=1$ for all places $v$ of $F$.

For fixed $F$, we can in constant (deterministic) time compute the set of even places of $F$.  We then compute the Hilbert symbol $(a,b)_v$ for $v$ real easily and for $v$ even by Algorithm \ref{evenhilbalg}.  

For the odd places, we first apply Lemma \ref{padicdet}, which implies that we need only check primes $\frakp \mid ab \Z_F$.  We compute $\frakg=a\Z_F+b\Z_F$ and then by small linear combinations we find $g \in \frakg^{-1}$ such that $g\frakg^{-1}$ is coprime to $a\Z_F$ and $b\Z_F$ and $(a+b)\Z_F$.  Now $\quat{a,b}{F} \cong \quat{a',b'}{F}$ where $a'=a+b$ and $b'=-abg^2$.  We claim that after repeating this eventually we will have $a$ and $b$ coprime.  Indeed, if $\ord_\frakp(a)=\ord_\frakp(b)$ then already $\ord_\frakp(-abg^2)=0$, and if $\ord_\frakp(a) > \ord_\frakp(b) > 0$, say, then $\ord_\frakp(-abg^2)=\ord_\frakp(a)-\ord_\frakp(b)$ and $\ord_\frakp(a+b)=\ord_\frakp(b)$, so then $\ord_\frakp(a)+\ord_\frakp(b) > \ord_\frakp(a)= \ord_\frakp(a')+\ord_\frakp(b')$, and since this is a sequence of nonnegative integers eventually either we will have either $\ord_\frakp(a)=0$ or $\ord_\frakp(b)=0$.

Then for any prime $\frakp \mid b\Z_F$, we have that $\frakp$ is ramified in $B$ if and only if $\frakp \mid \sqrad(b\Z_F)$ and $(a/\frakp)=-1$.  We can test this latter condition for all $\frakp \mid b\Z_F$ by calling the algorithm to solve (\textsf{QuadraticResiduosity}) by determining if $a$ is a quadratic residue modulo $\sqrad(b\Z_F)$.  We then repeat this step with $a,b$ interchanged, and we return \textsf{true} if and only if both of these quadratic residuosity tests return \textsf{true}.
\end{proof}

When $F=\Q$, in fact these problems are equivalent.

\begin{thm} \label{quadrecequiv}
Problem \textup{(\textsf{IsMatrixRing})} over $\Q$ is probablistic polynomial-time equivalent to Problem \textup{(\textsf{QuadraticResiduosity})} over $\Q$.
\end{thm}

\begin{rmk}
R\'onyai \cite{Ronyai1,RonyaiSimpleAlgebras} conditionally proves exactly Theorem \ref{quadrecequiv} (under the assumption of the Generalized Riemann Hypothesis).
\end{rmk}

Before proving this theorem, we derive one preliminary result.

\begin{lem} \label{qlb}
Let $a,b \in \Z_{>0}$ be such that $b$ is odd and $(a/b)=1$.  Let $\ell$ be an odd prime such that $\ell b \in (\Z/a\Z)^{\times 2}$ and $\displaystyle{\legen{a}{\ell}}=1$.  Then $\quat{a,\ell b}{\Q} \cong \M_2(\Q)$ if and only if $a$ is a square modulo $\sqrad(b)$.
\end{lem}

\begin{proof}
Again, we have $\quat{a,\ell b}{\Q} \cong \M_2(F)$ if and only if $(a,\ell b)_v=1$ for all places $v$ of $\Q$.  Since $a>0$, we know $(a,\ell b)_\infty=1$.  By hypothesis, for all odd $p \mid a$ we have $(\ell b/p)=1$ hence $(a,\ell b)_p=1$, and similarly $(a,\ell b)_\ell=1$.  Moreover, since $(a/b)=1$, the number of primes $p \mid \sqrad(b)$ such that $(a/p)=-1$ must be even, and since the quaternion algebra $\quat{a,\ell b}{\Q}$ is ramified at an even number of places, we conclude that $(a,\ell b)_2=1$.  Therefore $\quat{a,\ell b}{\Q} \cong \M_2(F)$ if and only if $(a,\ell b)_p=1$ for all $p \mid \sqrad(b)$ if and only if $a$ is a square modulo $\sqrad(b)$.
\end{proof}

The preceding lemma shows that the two problems in Theorem \ref{quadrecequiv} can be linked by finding a suitable prime $\ell$.  The conditions on $\ell$ are congruence conditions, so by the theorem on primes in arithmetic progression, such primes are abundant.  Explicitly, we rely on the specialization of a result from analytic number theory, stated by Adleman, Pomerance, and Rumely \cite[Proposition 8]{AdPomRum} and attributed to the proof of Linnik's theorem by Bombieri (using results of Gallagher and related to a result of Tatuzawa); see their paper for further discussion.  

\begin{lem} \label{largesieve}
There exist effectively computable (absolute) constants $x_0,\delta \in \R_{>0}$ such that whenever $x \geq x_0$, we have
\[ \Biggl| \sum_{\substack{\ell \leq x \\ \ell \equiv b \psmod{q}}} \log \ell - \frac{x}{\phi(q)}\Biggr| \leq \frac{x}{2 \phi(q)} \]
for all $q$ with $1 \leq q \leq x^{\delta}$ and all $b$ with $\gcd(b,q)=1$, except possibly for those $q$ which are multiples of a certain integer $q_0(x) > (\log x)^{3/2}$.
\end{lem}

\begin{proof}[Proof of Proposition \textup{\ref{quadrecequiv}}]
We must show that if we are able to solve (\textsf{IsMatrixRing}), then we can solve Problem (\textsf{QuadraticResiduosity}) in probabilistic polynomial time.  

Let $x = \max((4b)^{1/\delta}, x_0)$, with $x_0,\delta$ as in Lemma \ref{largesieve}.  Let $c$ be a random integer with $1 \leq c<b$.  We compute $q \equiv ac^2 \pmod{4b}$ with $1\leq q < 4b$ and $q \equiv 1 \pmod{4}$.  Then $q$ is a random element in $[1,4b] \cap \Z$ such that $aq \in (\Z/b\Z)^{\times 2}$ and $q \equiv 1 \pmod{4}$.  Let
\[ Q=\{1 \leq q<b: aq \in (\Z/b\Z)^{\times 2} \text{ and } q \equiv 1 \psmod{4}\}. \]

From Lemma \ref{largesieve}, we have $\sum_{\ell \leq x,\ \ell \equiv a \psmod{q}} \log \ell < x/(2 \phi(q))$ only if $q$ is divisible by $q_0(x) > (\log x)^{3/2}$; thus the set of such $q \in Q$ has cardinality at most $\# Q/(\log x)^{3/2}$.  Using partial summation (a standard argument which can be found in Davenport \cite[p.112]{Davenport}), it follows that a random $q \in Q$ has probability $1-1/(\log x)^{3/2}$ of satisfying 
\[ \pi(x;q,b)=\#\{ \ell \leq x: \text{$\ell$ prime},\ \ell \equiv b \psmod{q}\} < \frac{1}{2 \phi(q)}\frac{x}{\log x} \]
whenever $\gcd(b,q)=1$.  We then compute a random integer $1 \leq \ell < x$ with $\ell \equiv b \pmod{q}$ and test if $\ell$ is prime, which can be done in (deterministic) polynomial time \cite{AKS}.  Combining these, in probabilistic polynomial time, we may assume that $\ell$ indeed is prime.

We conclude by calling the algorithm to solve (\textsf{IsMatrixRing}) on $B=\quat{q,\ell b}{\Q}$.  We have
\[ \legen{q}{\ell}=\legen{\ell}{q}=\legen{b}{q}=\legen{q}{b}=\legen{a}{b}=1 \]
since $q \equiv 1\pmod{4}$, and $\ell b \equiv 1 \pmod{q}$.  So by Lemma \ref{qlb}, we have $B \cong \M_2(\Q)$ if and only if $q$ is a square modulo $\sqrad(b)$, which holds only if $a$ is a square modulo $\sqrad(b)$, as desired.
\end{proof}

We leave the natural generalization where $\Q$ is replaced by a number field $F$ as an open question.

\end{document}